\documentclass[11pt]{article}
\usepackage{amsmath,amssymb}
\usepackage{pstricks}
\usepackage{pst-node}
\usepackage{pst-poly}
\usepackage{multido}
\usepackage{tikz}
\newtheorem{result}{Theorem}

\newtheorem{define}{Definition}
\newtheorem{support}{Lemma}
\newtheorem{propo}{Proposition}

\newtheorem{note}{Remark}
\newtheorem{eg}{Example}

\newtheorem{observ}{Observation}

\newcommand{\qed}{%
\ifmmode 
\else \leavevmode\unskip\penalty9999 \hbox{}\nobreak\hfill \fi
\quad\hbox{\qedsymbol}}
\newcommand{\openbox}{\leavevmode \hbox to.77778em{%
\hfil\vrule
\vbox to.675em{\hrule width.6em\vfil\hrule}%
\vrule\hfil}}
\newcommand{\qedsymbol}{\openbox}
\newcommand{\showgrid}{}
\newcommand{\gridon}{\renewcommand{\showgrid}{\psset{subgriddiv=1,griddots=10,gridlabels=6pt}\psgrid}}
\gridon
\begin{document}

\begin{center} {\bf\LARGE On Edge-Partitioning of Complete Geometric Graphs into  Plane Trees}  \end{center}

\vskip8pt

\centerline{\large Hazim Michman  Trao$^{a}$, Gek L. Chia$^{b}$,  \  Niran Abbas Ali$^c$ \ and \   Adem Kilicman$^d$ }

\begin{center}
\itshape\small  $^{a, c, d}\/$Department of Mathematics, \\ Universiti Putra Malaysia, 43400 Serdang, Malaysia,  \\
 \vspace{1mm}
 $^{b}\/$Department of Mathematical and Actuarial Sciences, \\  Universiti Tunku Abdul Rahman, Sungai Long Campus,    Malaysia \\
\end{center}

\begin{abstract}
 In response to a well-known open question ``Does every complete geometric graph on $2n\/$ vertices have a partition of its edge set into $n\/$ plane spanning trees?" we provide an affirmative answer when the complete geometry graph is in  the regular wheel configuration. Also we present sufficient conditions for
 the complete geometric graph on $2n\/$ vertices to have a partition of its edge set into $n\/$ plane spanning trees (which are double stars, caterpillars or $                                                                                                                                         w\/$-caterpillars).

\end{abstract}

\section{Introduction}
A \emph{geometric graph} is a graph $G = (P,E)$ whose vertex set $P\/$ is a set of points in general position in the plane (that is, no three points lie on a common line) and whose edge set $E\/$ contains straight-line segments connecting the corresponding points.

\vspace{1mm} A {\em partition\/} of the edge set of a graph is a grouping of the edges into subgraphs so that every edge is in exactly one of the subgraphs.

\vspace{1mm} A well-known problem concerning  geometric graphs asks whether  every complete geometric graph on $2n\/$ vertices has a partition of its edge set into $n\/$ plane spanning trees?

\vspace{1mm}  In the case that the given graph is a convex complete graph on $2n\/$ vertices, an affirmative answer to the above question follows from a result of  Bernhart and Kanien \cite{bk:refer}.

\vspace{1mm}   In \cite{phrw:refer}, Bose et al. showed that if $T_1, T_2, \ldots, T_n\/$ is a partition of the edges of the complete convex {\red geometric} graphs $K_{2n}\/$, then   $T_1, T_2, \ldots, T_n\/$ are symmetric caterpillar that are pairwise isomorphic; and conversely for any symmetric convex caterpillar $T\/$ on $2n\/$ vertices, the edges of $K_{2n}\/$ can be partitioned into $n\/$ plane spanning convex copies of $T\/$ that are pairwise isomorphic. In addition, they also presented a sufficient condition for the edge set  of a complete geometric graph on $2n\/$ vertices  to be partitioned into $n\/$ plane
spanning double stars (that are pairwise graph-isomorphic).

\vspace{1mm} In Theorem \ref{the1}, we prove a similar result on complete geometric graph  $GW_{2n}\/$ in regular wheel configuration (with  $2n\/$ vertices). Here $T_1, T_2, \ldots, T_{n-1}\/$ are pairwise $w\/$-caterpillars and $T_n\/$ is a caterpillar.

 \vspace{1mm} Incidentally, in \cite{ahkv:refer}  Aichholzer et al.  proved that a complete geometric graph $GW_{2n}\/$ can be partitioned into $n\/$ plane spanning trees. The authors also noted that, in the case that $n\geq 3\/$, none of these trees can be paths.

 \vspace{1mm} In Theorem \ref{the2} we present a different sufficient condition for the edge set of a complete geometric graph on $2n\/$ vertices  to be partitioned into $n\/$ isomorphic copies of plane spanning trees (which are subdivision of double stars).

\vspace{1mm} In addition, we also   present a  sufficient condition for the edge set of a complete geometric graph on $2n\/$ vertices  to be partitioned into $n\/$ isomorphic copies of plane spanning symmetric caterpillars (see Theorem \ref{the3}). 

\vspace{1mm} Theorem \ref{the4}, we present a sufficient condition for the edge set of the complete geometric graph on $2n\/$ vertices to be partitioned into $n-1\/$ isomorphic copies of plane spanning symmetric $w\/$-caterpillars. This result contains Theorem \ref{the1} as a special case.

\vspace{1mm} We note in passing{\red ,} that Dor and Tarsi \cite{dt:refer} showed that the problem of deciding whether a given graph  can be partitioned into subgraphs isomorphic to a given graph is $NP\/$-complete.  Earlier,  Tutte \cite{wt:refer} and Nash-Williams \cite{nw:refer} independently obtained necessary and sufficient conditions for an abstract  graph to admit $k\/$ edge-disjoint spanning trees.
In \cite{sk:refer} Kundu  showed that any k-edge-connected graph contains at least $\lceil\frac{k-1}{2}\rceil\/$  edge-disjoint spanning trees.
\vspace{1mm} In the case that the trees are paths and stars, Priesler \cite{mt:refer} and Tarsi \cite{pt:refer} gave a characterization of partitioning the complete graph into paths and stars, respectively.


\section{Preliminaries}\label{sec-pre}

Let $P\/$ be a set of $2n\/$ points in the general position in the plane. A {\em complete\/} geometric graph $K_{2n}\/$ is a geometric graph on $P\/$ that has an edge joining every pair of points in $P\/$.

\vspace{1mm} A geometric graph is said to be \textit{plane} (or non-crossing) if its edges do not cross each other.

\vspace{1mm} A {\em tree\/} is a connected  graph with no cycle. Let $T\/$ denote a tree with $m\/$ vertices. It is easy to see that for any two vertices in $T\/$, there is a unique path joining these two vertices. If $T\/$ has only two vertices of degree $1\/$, then  $T\/$ is called an {\em $m\/$-path\/}.   Let $T'\/$ be  the tree obtained from $T\/$ by deleting all vertices of degree $1\/$. If $T'\/$ is a path, then $T\/$ is called a {\em caterpillar\/}. A {\em double star\/} is a caterpillar in which $T'\/$ is a $2\/$-path $vw\/$; if $d_{T}(v) = m+1\/$ and $d_{T}(w) = n+1\/$, then the double star is denoted by $S(m,n)\/$.

\vspace{1mm}  The \emph{convex hull} of $P\/$ is a smallest convex set that contains all points of $P\/$ and is denoted by $CH(P)\/$. A set  $P\/$ is  said to be in the \emph{convex position} in the plane if all points are on the boundary of $CH(P)\/$. A \emph{ convex geometric graph} is a  graph with the vertices in convex position.

\vspace{1mm} Suppose $G\/$ is a convex geometric graph on the set $P=\{0,1,...,2n-1\}\/$ where  the vertices are in anticlockwise order.  Let $G[i,j]\/$ denote the subgraph of $G\/$ induced by the vertices $\{v_i,v_{i+1},\dots ,v_j\}\/$ if $i<j\/$ and $\{v_i,v_{i+1},\dots ,v_0, v_1,v_2, \dots ,v_j\}\/$ if $i>j\/$ where all indices vertices are taken modulo $2n\/$.

\vspace{1mm}
\begin{support} \cite{phrw:refer} \label{lem2}
Let $v_iv_j\/$ be a non-boundary edge of a plane convex spanning tree $T\/$ such that $T[i,j]\/$ has exactly one boundary edge of $T\/$. Then exactly one of $v_iv_{j-1}\/$ and $v_jv_{i+1}\/$ is an edge of $T\/$.
\end{support}

Let $v_iv_j\/$ be a non-boundary edge of a plane convex spanning tree $T\/$ such that $T[i,j]\/$ has exactly one boundary edge of $T\/$. Then exactly one of $v_iv_{j-1}\/$ and $v_jv_{i+1}\/$ is an edge of $T\/$.

\vspace{1mm} An edge of a geometric graph $G\/$ on the boundary of the convex hull of $G\/$ is called a {\em boundary edge\/}. Two boundary edges are in {\em consecutive order\/} if they are incident the same vertex.

\vspace{1mm} Let $k\/$ be the number of pendant vertices of the tree  after removing all the pendant vertices of $T\/$. Garc\'{\i}a et al. \cite{ghhnt:refer} proved the following result:

\vspace{1mm}
\begin{support} \cite{ghhnt:refer} \label{lem1}
Let $T\/$ be a plane tree with at least two edges. In every convex drawing of $T\/$ there are at least $max\{2,k\}\/$ boundary edges. Moreover, if $T\/$ is not a star, then every plane convex drawing of $T\/$ has at least two non-consecutive boundary edges.
\end{support}

\vspace{1mm}
\section{Points in wheel configuration}\label{sec-wc}

A set $P\/$ of $m\/$, (is an even number) points is said to be in regular wheel configuration if $m-1\/$ of its points are regularly spaced on a circle $C\/$ with one point $x\/$ in the center of $C\/$. Note that those vertices in $C\/$ are the convex hull of $P\/$ and the anticlockwise ordering of the vertices of $C\/$ around the convex hull. An edge of the form $xv\/$ is called a radial edge; all other edges are called non-radial edges.

\vspace{1mm} In what follows, $GW_{2n}\/$ denotes the complete geometric graph on $2n\/$ vertices in wheel configuration.

\vspace{1mm}
\begin{observ} \label{ops1}
\vspace{1mm} Let $P\/$ be a set of $2n\/$ points in regular wheel configuration in the plane where $n \geq 2\/$. Suppose $uv\/$ is a non-radial edge of a plane spanning tree $T\/$. Then $uv\/$ separates $P-\{u, v\}\/$ into two parts $A\/$ and $A(x)\/$ with (i) $A(x)\/$ containing the center $x\/$, (ii) $|A| < |A(x)|\/$ and (iii)  $|A| \leq n-2\/$ and $|A(x)| \geq n\/$.
\end{observ}

\vspace{1mm}
\begin{propo} \label{pro1}
Let $P\/$ be a set of $2n\/$ points in regular wheel configuration in the plane where $n \geq 2\/$. Suppose  a plane spanning tree $T\/$ on $P\/$ has no boundary edges. Then all edges in $T\/$ are radial edges.
\end{propo}

\vspace{1mm}  \noindent
{\bf Proof:} Let $v_0, v_1, \ldots, v_{2n-2}\/$ be the vertices of $C\/$ in cyclic order. Suppose on the contrary that $T\/$ has a non-radial edge $v_iv_j\/$. Then by Lemma \ref{lem2},  $T[i,j]\/$ (which is in convex position)  has at least two boundary edges one of which is different from  $v_iv_j\/$, a contradiction. \qed

\vspace{2mm} An immediate consequence of the above result is the following lemma.

\vspace{1mm}
\begin{support} \label{lem3}
Let $P\/$ be a set of $2n\/$ points in regular wheel configuration in the plane where $n \geq 2\/$. Suppose $T_1,T_2,...,T_n\/$ is a partition of edges of $GW_{2n}$ on $P\/$ into plane spanning trees. Then each $T_i$, $i=1,2,...,n\/$ has at least one boundary edge.
\end{support}

\vspace{1mm}  \noindent
{\bf Proof:} Suppose $x\/$ is the center of $P\/$ and $v_0, v, v_1,  \cdots, v_{2n-2}\/$ are the vertices of the circle $C\/$. Assume $T_1,T_2,...,T_n\/$ is a partition of $GW_{2n}$ on $P\/$ into plane spanning trees.\

\vspace{1mm}
\vspace{1mm} Assume to the contrary that $T_k\/$ has no boundary edge for some $1\leq k\leq n\/$. By Proposition \ref{pro1} all radial edges $xv_i\/$ for $i=0,1,...,2n-1\/$ are in $T_k\/$, a contradiction since each of $\{T_1,T_2,...,T_n\}\/$ are spanning trees. \qed

\vspace{1mm}
\begin{propo} \label{pro2}
Let $P\/$ be a set of $2n\/$ points in regular wheel configuration in the plane where $n \geq 2\/$. Suppose  a plane spanning tree $T\/$ has only one boundary edge. Then $T\/$ has at least $n\/$ consecutive radial edges.
\end{propo}

\vspace{1mm}  \noindent
{\bf Proof:} Suppose $x\/$ is the center of $P\/$ and $v_0, v, v_1,  \cdots, v_{2n-2}\/$ are the vertices of the circle $C\/$.

\vspace{1mm} By Observation \ref{ops1}, $T$ has a non-radial edge $uv$ that separates $P-\{u, v\}$ into $A$ and $A(x)$ with $|A|<|A(x)|$. Choose $uv$ such that $A$ is maximal with respect to this property.

\vspace{1mm} Relabel the vertices of $C\/$  such that $u=v_0\/$ and $v= v_i \/$. The choice on $uv\/$ implies that $v_0v_j \not \in E(T)\/$ for any $j\/$ such that $j >i\/$ where $i \leq n-1\/$.

\vspace{1mm} Since $T\/$ has only one boundary edge, by Lemma \ref{lem1}, the convex subtree $T'\/$ induced by the set of vertices $\{v_0, v_1, \ldots, v_i\}\/$ has only two boundary edges edges, one of which is $v_0v_i\/$.

\vspace{1mm} Clearly, either $xv_0 \in E(T)\/$ or else $xv_i \in E(T)\/$.   Consider the plane subtree $T"\/$ of $T\/$ induced by the set of vertices $\{ w, v_{i+1}, v_{i+2}, \ldots, v_{2n-2}, x\}\/$ where $w \in \{v_0, v_i\}$. By Proposition \ref{pro1}, all edges of $T"\/$ are radial edges.

\vspace{1mm} It follows that $T = T' \cup T"\/$ as at least $n\/$ radial edges in consecutive order.   \qed

\vspace{2mm} Immediate from Propositions \ref{pro1} and \ref{pro2} is the following result.

\vspace{1mm}
\begin{propo} \label{pro3}
Let $P\/$ be a set of $2n\/$ points in regular wheel configuration in the plane where $n \geq 2\/$. Suppose $T\/$ is a plane spanning tree has $k\/$ radial edges. Then $T\/$ has at least two boundary edges if $k<n\/$.
\end{propo}

\vspace{1mm}
\begin{support} \label{lem4}
Let $P\/$ be a set of $2n\/$ points in regular wheel configuration in the plane where $n \geq 2\/$. Suppose $T_1,T_2,...,T_n\/$ is a partition of $GW_{2n}$ on $P\/$ into plane spanning trees. Then
there is a unique plane spanning tree  $T\in \{T_1,T_2,...,T_n\}\/$ such that (i) $T\/$ has only one boundary edge with only $n\/$ radial edges  in consecutive order, and (ii) each tree in $\{T_1,T_2,...,T_n\}-T\/$ has exactly two boundary edges and one radial edge.
\end{support}

\vspace{1mm}  \noindent
{\bf Proof:}  By Lemma \ref{lem3} each $T_i\/$ has at least one boundary edge.

\vspace{1mm} Since $GW_{2n}\/$ has only $2n-1\/$ boundary edges, it is not possible that each $T_i\/$ has $2\/$  boundary edges.  Hence some plane spanning tree $T\/$ from $\{T_1, T_2, \ldots, T_n\}\/$ has only one boundary edge.

\vspace{1mm} By Proposition \ref{pro2}, $T\/$ has at least $n \/$ radial edges  (in consecutive order). Since $GW_{2n}\/$ has only $2n-1\/$ radial edges, and each plane spanning tree from $\{T_1, T_2, \ldots, T_n\} - T\/$ has at least one radial edge, it follows that $T\/$ has exactly $n\/$ radial edges. This also implies that each plane spanning tree from $\{T_1, T_2, \ldots, T_n\} - T\/$ has only one radial edge and two boundary edges (by Proposition \ref{pro3}). This completes the proof. \qed

\vspace{1mm}
Suppose $T\/$ is a tree and $u, v\/$ are two vertices in $T\/$. We let $P(u,v)\/$ denote the unique path from $u\/$ to $v\/$ in $T\/$.

\vspace{1mm}
\begin{support} \label{lem5}
Let $P\/$ be a set of $2n\/$ points in regular wheel configuration in the plane where $n \geq 2\/$. Suppose $T_1,T_2,..., T_n\/$ is a partition of $GW_{2n}$ on $P\/$ into plane spanning trees where $T_i\/$ has two boundary edges, $ i=1, \ldots, n-1\/$. Suppose $xv_i \in E(T_i)\/$.   Then
for each $i=1,2, \dots , n-1\/$, $xv_iv_{i+n}v_{i+1}\/$ (or equivalently $xv_iv_{i+n-1}v_{i-1}\/$) is a path in $T_i\/$.
\end{support}

\vspace{1mm}  \noindent
{\bf Proof:}  By Lemma \ref{lem4}, $T_n\/$ has $n\/$ consecutive radial edges, say given by  $v_ix\/$, $i= n, n+1, \ldots, 2n-2,  0\/$.

\vspace{1mm} By Lemma \ref{lem4}, $T_i\/$ has only one radial edge for each $i =1, \ldots, n-1\/$. Without loss of generality assume that $xv_i\/$ is a radial edge of $T_i\/$, $i=1, \ldots, n-1\/$.

\vspace{1mm} Consider the unique path $P(v_{1+n}, v_1) \/$ in $T_1\/$ from $v_{1+n}\/$ to $v_1\/$. Suppose $P(v_{1+n}, v_1)\/$ has length $r\/$. Since $T_1\/$ is a plane tree, it follows that $T_1\/$ has $r\/$ boundary edges that lie on the segment from $v_{1+n}\/$ to $v_1\/$ (in anti-clockwise direction).

\vspace{1mm} Likewise, if the unique path $P(v_{1+n}, v_2)\/$ in $T_1\/$ from $v_2\/$ to $v_{1+n}\/$  has length $s\/$, then $T_1\/$ has $s\/$ boundary edges that lie in the segment from $v_2\/$ to $v_{1+n}\/$  (in anti-clockwise direction). Since $T_1\/$ has only $2\/$ boundary edges, it follows that $r = 1= s\/$. Hence $xv_1v_{1+n}v_2\/$ is a path in $T_1\/$.

\vspace{1mm} By applying similar argument to each of $T_2, \ldots, T_{n-1}\/$, we have the conclusion that $xv_iv_{i+n}v_{i+1}\/$ is a path in $T_i\/$ for each $i=1, 2, \ldots, n-1\/$.   \qed

\vspace{2mm}
In \cite{phrw:refer}, the authors call a tree $T\/$ is symmetric if there exists an edge $vw\/$ of $T\/$ such that if $A\/$ and $B\/$ are the components of $T-\{vw\}\/$ with $v\in A\/$ and $w\in B\/$, then there exists a graph-isomorphism between $A\/$ and $B\/$ that maps $v\/$ to $w\/$.

\vspace{1mm}
\begin{define} \label{def1}
Let $T$ be a tree with $2n\/$ vertices containing a $4$-path $u_1u_2u_3u_4\/$ such that the resulting graph obtained from $T$ by deleting the edges of the $4$-path consists of two isolated vertices $u_1,u_i\/$ for some $i\in \{2,4\}\/$ together with two non-trivial components $A,B\/$ where $u_3\in A\/$ and $u_{6-i}\in B\/$. Then $T\/$ is said to be $P_4\/$-symmetric if there is a graph-isomorphism   $\varphi\/$ between $A\/$ and $B\/$ on $T - \{u_1u_2, u_2u_3, u_3u_4\}\/$  such that $\varphi(u_3)=u_{6-i}\/$. When $n \geq 3\/$ and  $i=4\/$ (respectively $i=2\/$), the $P_4\/$-symmetric tree $T\/$ said to be of Type-1 (respectively Type-2).
\end{define}

\begin{figure}[htb]
\centering
\resizebox{11cm}{!}{
\begin{minipage}{.45\textwidth}
\begin{tikzpicture}

       \coordinate (u1) at (0,0);\filldraw[black] (u1) circle(3.5pt);\node[left] at (u1) {${u_1}$};
       \coordinate (v0) at (-0.55,-2.45);\filldraw[black] (v0) circle(3.5pt);
       \coordinate (v1) at (1.56,-1.96);\filldraw[black] (v1) circle(3.5pt);
       \coordinate (u2) at (2.5,0);\filldraw[black] (u2) circle(3.5pt);\node[below right] at (u2) {${u_2}$};
       \coordinate (u4) at (1.56,1.96);\filldraw[black] (u4) circle(3.5pt);\node[above right] at (u4) {${u_4}$};
       \coordinate (v4) at (-0.55,2.45);\filldraw[black] (v4) circle(3.5pt);
       \coordinate (v5) at (-2.25,1.08);\filldraw[black] (v5) circle(3.5pt);
       \coordinate (u3) at (-2.25,-1.08);\filldraw[black] (u3) circle(3.5pt);\node[below left] at (u3) {${u_3}$};
       \draw [line width=1,black](v0) -- (v1);\draw [line width=1,black](v0) -- (u2);\draw [line width=1,red](u2) -- (u1);
       \draw [line width=1,red](u2) -- (u3);\draw [line width=1 ,red](u3) -- (u4);\draw [line width=1,black](u3) -- (v4);
       \draw [line width=1,black](v4) -- (v5);

     \filldraw[red] (u1) circle(1.5pt);\filldraw[white] (v0) circle(1.5pt);\filldraw[white] (v1) circle(1.5pt);\filldraw[red] (u2) circle(1.5pt);
     \filldraw[red] (u4) circle(1.5pt);\filldraw[white] (v4) circle(1.5pt);\filldraw[white] (v5) circle(1.5pt);\filldraw[red] (u3) circle(1.5pt);

\end{tikzpicture}

\centering
$T_1\/$
\end{minipage}
\begin{minipage}{.45\textwidth}
\begin{tikzpicture}

       \coordinate (u1) at (0,0);\filldraw[black] (u1) circle(3.5pt);\node[left] at (u1) {${u_1}$};
       \coordinate (v0) at (-0.55,-2.45);\filldraw[black] (v0) circle(3.5pt);
       \coordinate (v1) at (1.56,-1.96);\filldraw[black] (v1) circle(3.5pt);
       \coordinate (u2) at (2.5,0);\filldraw[black] (u2) circle(3.5pt);\node[below right] at (u2) {${u_2}$};
       \coordinate (u4) at (1.56,1.96);\filldraw[black] (u4) circle(3.5pt);\node[above right] at (u4) {${u_4}$};
       \coordinate (v4) at (-0.55,2.45);\filldraw[black] (v4) circle(3.5pt);
       \coordinate (v5) at (-2.25,1.08);\filldraw[black] (v5) circle(3.5pt);
       \coordinate (u3) at (-2.25,-1.08);\filldraw[black] (u3) circle(3.5pt);\node[below left] at (u3) {${u_3}$};

       \draw [line width=1,black](v0) -- (v1);\draw [line width=1,black](v1) -- (u3);\draw [line width=1,red](u2) -- (u1);
       \draw [line width=1,red](u2) -- (u3);\draw [line width=1,red](u3) -- (u4);\draw [line width=1,black](u4) -- (v5);
       \draw [line width=1,black](v4) -- (v5);

     \filldraw[red] (u1) circle(2pt);\filldraw[white] (v0) circle(2pt);\filldraw[white] (v1) circle(2pt);\filldraw[red] (u2) circle(2pt);
     \filldraw[red] (u4) circle(2pt);\filldraw[white] (v4) circle(2pt);\filldraw[white] (v5) circle(2pt);\filldraw[red] (u3) circle(2pt);

\end{tikzpicture}

\centering
$T_2\/$
\end{minipage}
}

\vspace{5mm}

\caption{\small { Two $P_4\/$-symmetric trees $T_1$ and $T_2$. }}
 \label{reg}
\end{figure}
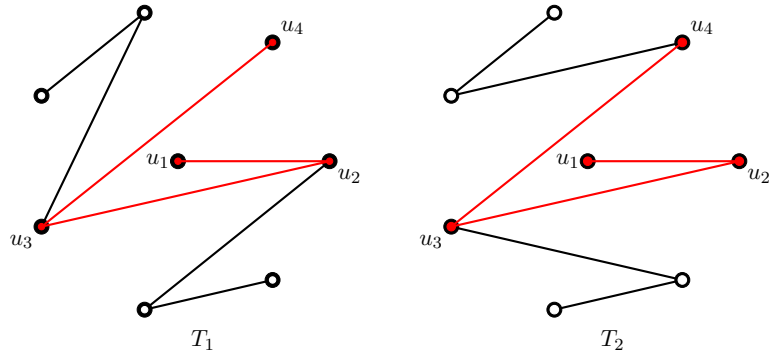

\vspace{1mm}
\begin{eg}\label{eg1}
 Figure \ref{reg} depicts two $P_4\/$-symmetric trees $T_1\/$ and $T_2\/$. Here both $T_1\/$ and $T_2\/$ contains a $4\/$-path $u_1u_2u_3u_4\/$. Note also that  in $T_1 -\{u_1u_2,u_2u_3, u_3u_4\}$, we have $\varphi(u_3)=u_{2}\/$,  and in  $T_2-\{u_1u_2,u_2u_3, u_3u_4\}$ we have  $\varphi(u_3)=u_{4}\/$.  Here $T_1\/$ is Type-1 and $T_2\/$ is Type-2.
\end{eg}

\vspace{1mm}

\vspace{1mm}
\begin{define} \label{def2}
Let $T\/$ be a $P_4\/$-symmetric tree with a $4\/$-path  $u_1u_2u_3u_4\/$. Then we call $T\/$ a $w\/$-caterpillar  if either  $T\/$ is a caterpillar or else $T -u_1\/$ is a caterpillar.
\end{define}

\vspace{1mm} The two $P_4\/$-symmetric trees depicted in Figure \ref{reg} are both $w\/$-caterpillars.

\vspace{1mm}
\begin{result} \label{the1}
Let $P\/$ be a set of $2n$ points in regular wheel configuration in the plane where $n \geq 2\/$. Suppose $\/T_1,T_2, \ldots,T_n$ is a partition of the edges of the complete graph $GW_{2n}$ on $P\/$ into plane spanning trees. Then  $\/n-1\/$ of these trees $T_1, \ldots, T_{n-1}\/$  are $w\/$-caterpillars and  are pairwise isomorphic, and  $T_n\/$ is a  caterpillar with only one boundary edge.
 Conversely, for any  $w\/$-caterpillar $T\/$,   the edges of the complete graph $GW_{2n}\/$ can be partitioned into $n-1\/$ copies of $T$ together with a plane spanning caterpillar $T'\/$ with only one boundary edge.
\end{result}

\vspace{1mm}  \noindent
{\bf Proof:} The result is clearly true when $n=2\/$. Hence we assume that $n \geq 3\/$.

\vspace{1mm} By Lemma \ref{lem4}, $T_n\/$ has $n\/$ consecutive radial edges, say given by  $v_ix\/$, $i= n, n+1, \ldots, 2n-2, 0\/$. Moreover each $T_i\/$ has exactly two boundary edges, $i=1, 2, \ldots, n-1\/$.

\vspace{1mm} Suppose $xv_i \in E(T_i)\/$, $i =1, \ldots, n-1\/$.
 By Lemma \ref{lem5}, $T_i\/$ has a 4-path $xv_iv_{i+n}v_{i+1}\/$.

 \vspace{1mm} Note that the edge $v_1v_n\/$ crosses an edge in the  $4\/$-path $xv_iv_{i+n}v_{i+1}\/$ for each $i =1, \ldots, n-1\/$. This implies that  $v_1v_n \in E(T_n)\/$.

\vspace{1mm} Since each $T_i\/$ has exactly two boundary edges, and since $v_{i+n}v_i \in E(T_i)\/$ and $v_{i+1}v_{i+n} \in E(T_i)\/$,     it follows that $T_i\/$ has one boundary edge that lie on the segment from $v_{i+n}\/$ to $v_i\/$ and one boundary edge on the segment from $v_{i+1}\/$ to $v_{i+n}\/$ (both in anti-clockwise direction).

\vspace{1mm} To extend each of these $4\/$-paths to the required $w$-caterpillar, we first consider the plane tree $T_1\/$.

\vspace{1mm} Either (i) there is a unique path $P(v_{n+2}, v_1) \/$ in $T_1\/$ from $v_{n+2}\/$ to $v_1\/$ or (ii) there is a unique path $P(v_{n+1}, v_0) \/$ in $T_1\/$ from $v_{n+1}\/$ to $v_0\/$. Here the operations on the subscripts are reduced modulo $2n-1\/$.

\vspace{1mm} (i) Suppose   $P(v_{n+2}, v_1) \/$ is a subpath in $T_1\/$. Then   $v_1v_{n+2}\in E(T_1)\/$ (since $T_1\/$ is a plane tree).

\vspace{1mm} This implies that  $P(v_{n+2}, v_1) \/$ is a not a subpath of the plane tree $T_2\/$. Hence $P(v_{n+3}, v_2) \/$ is a subpath of $T_2\/$.  Consequently,  $v_2v_{n+3} \in E(T_2)\/$.

\vspace{1mm} Continue with the same argument, we have $v_iv_{n+i+1}\in E(T_i)\/$ for each $i =1, \ldots, n\/$.

\begin{figure}[htb]
\centering
\resizebox{12cm}{!}{
\begin{minipage}{.45\textwidth}
\begin{tikzpicture}

       \coordinate (x) at (0,0);\node at (x) {\textbullet};
       \coordinate (v0) at (0.78,-2.38);\filldraw[black] (v0) circle(2pt);\node[below] at (v0) {${v_0}$};
       \coordinate (v1) at (2.5,0);\filldraw[black] (v1) circle(2pt);\node[right] at (v1) {${v_1}$};
       \coordinate (v2) at (0.78,2.38);\filldraw[black] (v2) circle(2pt);\node[above] at (v2) {${v_2}$};
       \coordinate (v3) at (-2.02,1.48);\filldraw[black] (v3) circle(2pt);\node[above left] at (v3) {${v_3}$};
       \coordinate (v4) at (-2.02,-1.48);\filldraw[black] (v4) circle(2pt);\node[below left] at (v4) {${v_4}$};
    \draw [line width=1,red](v2) -- (x); \draw [line width=1,red](v2) -- (v0); \draw [line width=1,red](v0) -- (v3);
    \draw [line width=1,blue](v1) -- (x); \draw [line width=1,blue](v1) -- (v4); \draw [line width=1,blue](v4) -- (v2);
    \draw [line width=1,black](v3) -- (x);\draw [line width=1,black](v4) -- (x);\draw [line width=1,black](v0) -- (x);
     \draw [line width=1,black](v1) -- (v3);
    \draw [line width=1,blue](v0) -- (v4);\draw [line width=1,red](v3) -- (v4);\draw [line width=1,blue](v2) -- (v3);
    \draw [line width=1,black](v1) -- (v2);\draw [line width=1,red](v0) -- (v1);

       \filldraw[black] (x) circle(3pt);\filldraw[black] (v0) circle(3pt);\filldraw[black] (v1) circle(3pt);\filldraw[black] (v2) circle(3pt);
       \filldraw[black] (v3) circle(3pt);\filldraw[black] (v4) circle(3pt);

       \filldraw[white] (v0) circle(1.5pt);\filldraw[white] (v1) circle(1.5pt);\filldraw[white] (v2) circle(1.5pt);
       \filldraw[white] (v3) circle(1.5pt);\filldraw[white] (v4) circle(1.5pt);\node[below right] at (x) {${x}$};

\end{tikzpicture}
\centering

(i)

\end{minipage}
\begin{minipage}{.45\textwidth}
\begin{tikzpicture}

       \coordinate (x) at (0,0);\node at (x) {\textbullet};
       \coordinate (v0) at (0.78,-2.38);\filldraw[black] (v0) circle(2pt);\node[below] at (v0) {${v_0}$};
       \coordinate (v1) at (2.5,0);\filldraw[black] (v1) circle(2pt);\node[right] at (v1) {${v_1}$};
       \coordinate (v2) at (0.78,2.38);\filldraw[black] (v2) circle(2pt);\node[above] at (v2) {${v_2}$};
       \coordinate (v3) at (-2.02,1.48);\filldraw[black] (v3) circle(2pt);\node[above left] at (v3) {${v_3}$};
       \coordinate (v4) at (-2.02,-1.48);\filldraw[black] (v4) circle(2pt);\node[below left] at (v4) {${v_4}$};
    \draw [line width=1,red](v2) -- (x); \draw [line width=1,red](v2) -- (v0); \draw [line width=1,red](v0) -- (v3);
    \draw [line width=1,blue](v1) -- (x); \draw [line width=1,blue](v1) -- (v4); \draw [line width=1,blue](v4) -- (v2);
    \draw [line width=1,black](v3) -- (x);\draw [line width=1,black](v4) -- (x);\draw [line width=1,black](v0) -- (x);
     \draw [line width=1,black](v1) -- (v3);
     \draw [line width=1,blue](v1) -- (v0);\draw [line width=1,red](v1) -- (v2);\draw [line width=1,black](v2) -- (v3);
     \draw [line width=1,blue](v4) -- (v3);\draw [line width=1,red](v0) -- (v4);

       \filldraw[black] (x) circle(3pt);\filldraw[black] (v0) circle(3pt);\filldraw[black] (v1) circle(3pt);\filldraw[black] (v2) circle(3pt);
       \filldraw[black] (v3) circle(3pt);\filldraw[black] (v4) circle(3pt);

       \filldraw[white] (v0) circle(1.5pt);\filldraw[white] (v1) circle(1.5pt);\filldraw[white] (v2) circle(1.5pt);
       \filldraw[white] (v3) circle(1.5pt);\filldraw[white] (v4) circle(1.5pt);\node[below right] at (x) {${x}$};

\end{tikzpicture}
\centering

(ii)

\end{minipage}
\vspace{5mm}
}
\caption{An illustration of case (i) and case (ii) when $n=3\/$. }  \label{fig-2-th1}

\end{figure}
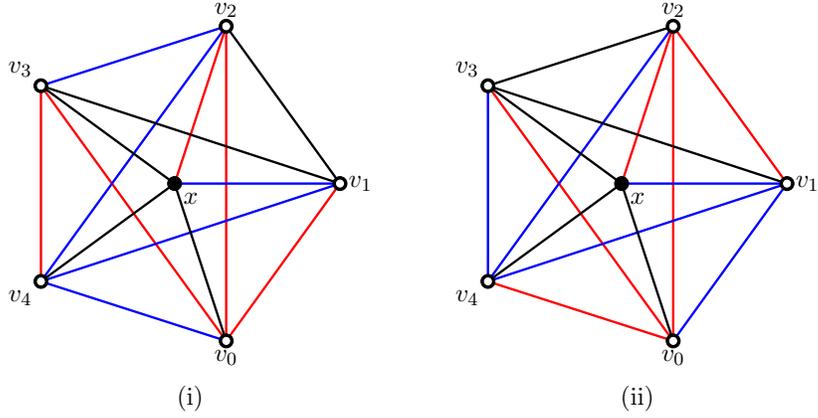

\vspace{1mm} Since $v_nv_{2n+1} = v_n v_2\/$ is an edge in $T_n\/$, it follows that  $v_{n+1}v_3\/$ is an edge in $T_1\/$. As such,   by applying the same argument as before, we conclude that  $v_{n+i}v_{i+2}\in E(T_{i})\/$ for each $i=1, 2, \ldots, n-1\/$. If $n=3\/$, we see that $T_1, T_2\/$ are Type-1 $P_4\/$-symmetric plane trees.  Clearly $T_1, T_2\/$ are caterpillars (see Figure \ref{fig-2-th1}).

\vspace{1mm} On the other hand, if (ii) there is a unique path $P(v_{n+1}, v_0) \/$ in $T_1\/$ from $v_{n+1}\/$ to $v_0\/$, then
$v_0v_{n+1}\in E(T_1)\/$ (since $T_1\/$ is a plane tree). But this implies that  $P(v_{n+1}, v_0) \/$ is a not a subpath of the plane tree $T_{n-1}\/$. Hence $P(v_{n}, v_{2n-2}) \/$ is a subpath of $T_{n-1}\/$.  Consequently,  $v_nv_{2n-2} \in E(T_{n-1})\/$.

\vspace{1mm} Continue with the same argument, we have $v_{i+1}v_{n+i-1}\in E(T_i)\/$ for each $i =n-1,n-2 \ldots,  0\/$. Here $T_0 = T_n\/$.

\vspace{1mm} Since $v_1 v_{n-1}\/$ is an edge in $T_n\/$, it follows that  $v_{n-2}v_0\/$  is an edge in $T_{n-1}\/$. As such,  by applying the same argument as before, we conclude that  $v_{i-1}v_{n+i}\in E(T_{i})\/$ for each $i =1, 2, \ldots, n-1\/$.  Again when $n=3\/$, we see that $T_1, T_2\/$ are Type-2 $P4\/$-symmetric plane trees. Clearly $T_1, T_2\/$ are caterpillars (see Figure \ref{fig-2-th1}).

\pagebreak
\vspace{1mm} Hence assume that $n\geq 4\/$. There are two cases to consider.

\vspace{1mm} {\em Case (a) $v_1v_{n+2} \in E(T_1)\/$  }

\vspace{1mm} By the argument in (i), $v_3v_{n+1} \in E(T_1)\/$.

\vspace{1mm} (a1) If $v_1v_{n+3} \in E(T_1)\/$, then $v_{n+1}v_4 \in E(T_1)\/$. Moreover $v_iv_{n+i+2}, v_{n+i}v_{i+3} \in E(T_i)\/$, $i=1, 2, \ldots, n-1\/$.

\vspace{1mm} (a2) If $v_0v_{n+2} \in E(T_1)\/$, then $v_{n}v_3 \in E(T_1)\/$. Moreover $v_{i-1}v_{n+i+1}, v_{n+i-1}v_{i+2} \in E(T_i)\/$, $i=1, 2, \ldots, n-1\/$.

\vspace{1mm} Continue with this argument, we see that each of the plane trees $T_1, T_2, \ldots, T_{n-1}\/$ are $P_4\/$-symmetric of Type-1 and that they are pairwise isomorphic. Here the $4\/$-path in $T_i\/$ is $xv_iv_{n+i}v_{i+1}\/$ and  $\varphi(v_i) = v_{n+i}\/$. It is easy to see that $T_i\/$ is a caterpillar.

\begin{figure}[htb]
\centering
\resizebox{12cm}{!}{
\begin{minipage}{.45\textwidth}
\begin{tikzpicture}

       \coordinate (x) at (0,0);\node at (x) {\textbullet};
       \coordinate (v0) at (1.56,-1.96);\filldraw[black] (v0) circle(2pt);\node[below right] at (v0) {${v_0}$};
       \coordinate (v1) at (2.5,0);\filldraw[black] (v1) circle(2pt);\node[right] at (v1) {${v_1}$};
       \coordinate (v2) at (1.56,1.96);\filldraw[black] (v2) circle(2pt);\node[above right] at (v2) {${v_2}$};
       \coordinate (v3) at (-0.56,2.45);\filldraw[black] (v3) circle(2pt);\node[above] at (v3) {${v_3}$};
       \coordinate (v4) at (-2.25,1.09);\filldraw[black] (v4) circle(2pt);\node[above left] at (v4) {${v_4}$};
       \coordinate (v5) at (-2.25,-1.09);\filldraw[black] (v5) circle(2pt);\node[below left] at (v5) {${v_5}$};
       \coordinate (v6) at (-0.56,-2.45);\filldraw[black] (v6) circle(2pt);\node[below] at (v6) {${v_6}$};

       \draw [line width=1,red](v1) -- (x);\draw [line width=1,red](v1) -- (v5);\draw [line width=1,red](v5) -- (v2);
       \draw [line width=1,blue](v2) -- (x);\draw [line width=1,blue](v2) -- (v6);\draw [line width=1,blue](v6) -- (v3);
       \draw [line width=1,green](v3) -- (x);\draw [line width=1,green](v3) -- (v0);\draw [line width=1,green](v0) -- (v4);
       \draw [line width=1,black](v4) -- (x);\draw [line width=1,black](v5) -- (x);\draw [line width=1,black](v6) -- (x);
       \draw [line width=1,black](v0) -- (x);\draw [line width=1,black](v1) -- (v4);
       \draw [line width=1,red](v1) -- (v6); \draw [line width=1,red](v5) -- (v3);
       \draw [line width=1,blue](v2) -- (v0); \draw [line width=1,blue](v6) -- (v4);
       \draw [line width=1,green](v3) -- (v1);\draw [line width=1,green](v0) -- (v5);
       \draw [line width=1,black](v2) -- (v4);

        \draw [line width=1,red](v1) -- (v0);\draw [line width=1,red](v5) -- (v4);
       \draw [line width=1,blue](v2) -- (v1);\draw [line width=1,blue](v6) -- (v5);
       \draw [line width=1,green](v3) -- (v2);\draw [line width=1,green](v0) -- (v6);
       \draw [line width=1,black](v4) -- (v3);

       \filldraw[black] (v0) circle(3pt);\filldraw[black] (v1) circle(3pt);\filldraw[black] (v2) circle(3pt);
       \filldraw[black] (v3) circle(3pt);\filldraw[black] (v4) circle(3pt);\filldraw[black] (v5) circle(3pt);\filldraw[black] (v5) circle(3pt);

       \filldraw[white] (v0) circle(1.5pt);\filldraw[white] (v1) circle(1.5pt);\filldraw[white] (v2) circle(1.5pt);
       \filldraw[white] (v3) circle(1.5pt);\filldraw[white] (v4) circle(1.5pt);\filldraw[white] (v5) circle(1.5pt);\filldraw[white] (v6) circle(1.5pt);
       \filldraw[black] (x) circle(3pt);\node[below right] at (x) {${x}$};

\end{tikzpicture}
\centering

(a1)

\end{minipage}
\begin{minipage}{.45\textwidth}
\begin{tikzpicture}

       \coordinate (x) at (0,0);\node at (x) {\textbullet};
       \coordinate (v0) at (1.56,-1.96);\filldraw[black] (v0) circle(2pt);\node[below right] at (v0) {${v_0}$};
       \coordinate (v1) at (2.5,0);\filldraw[black] (v1) circle(2pt);\node[right] at (v1) {${v_1}$};
       \coordinate (v2) at (1.56,1.96);\filldraw[black] (v2) circle(2pt);\node[above right] at (v2) {${v_2}$};
       \coordinate (v3) at (-0.56,2.45);\filldraw[black] (v3) circle(2pt);\node[above] at (v3) {${v_3}$};
       \coordinate (v4) at (-2.25,1.09);\filldraw[black] (v4) circle(2pt);\node[above left] at (v4) {${v_4}$};
       \coordinate (v5) at (-2.25,-1.09);\filldraw[black] (v5) circle(2pt);\node[below left] at (v5) {${v_5}$};
       \coordinate (v6) at (-0.56,-2.45);\filldraw[black] (v6) circle(2pt);\node[below] at (v6) {${v_6}$};

       \draw [line width=1,red](v1) -- (x);\draw [line width=1,red](v1) -- (v5);\draw [line width=1,red](v5) -- (v2);
       \draw [line width=1,blue](v2) -- (x);\draw [line width=1,blue](v2) -- (v6);\draw [line width=1,blue](v6) -- (v3);
       \draw [line width=1,green](v3) -- (x);\draw [line width=1,green](v3) -- (v0);\draw [line width=1,green](v0) -- (v4);
       \draw [line width=1,black](v4) -- (x);\draw [line width=1,black](v5) -- (x);\draw [line width=1,black](v6) -- (x);
       \draw [line width=1,black](v0) -- (x);\draw [line width=1,black](v1) -- (v4);
       \draw [line width=1,red](v1) -- (v6); \draw [line width=1,red](v5) -- (v3);
       \draw [line width=1,blue](v2) -- (v0); \draw [line width=1,blue](v6) -- (v4);
       \draw [line width=1,green](v3) -- (v1);\draw [line width=1,green](v0) -- (v5);
       \draw [line width=1,black](v2) -- (v4);
       \draw [line width=1,red](v0) -- (v6);\draw [line width=1,red](v3) -- (v4);
       \draw [line width=1,blue](v0) -- (v1);\draw [line width=1,blue](v4) -- (v5);
       \draw [line width=1,green](v1) -- (v2);\draw [line width=1,green](v5) -- (v6);
       \draw [line width=1,black](v2) -- (v3);

       \filldraw[black] (v0) circle(3pt);\filldraw[black] (v1) circle(3pt);\filldraw[black] (v2) circle(3pt);
       \filldraw[black] (v3) circle(3pt);\filldraw[black] (v4) circle(3pt);\filldraw[black] (v5) circle(3pt);\filldraw[black] (v5) circle(3pt);

       \filldraw[white] (v0) circle(1.5pt);\filldraw[white] (v1) circle(1.5pt);\filldraw[white] (v2) circle(1.5pt);
       \filldraw[white] (v3) circle(1.5pt);\filldraw[white] (v4) circle(1.5pt);\filldraw[white] (v5) circle(1.5pt);\filldraw[white] (v6) circle(1.5pt);
       \filldraw[black] (x) circle(3pt);\node[below right] at (x) {${x}$};

\end{tikzpicture}
\centering

(a2)

\end{minipage}
}

\vspace{5mm}

\caption{\small { An illustration of Case (a). }}
 \label{fig-3-th1}
\end{figure}
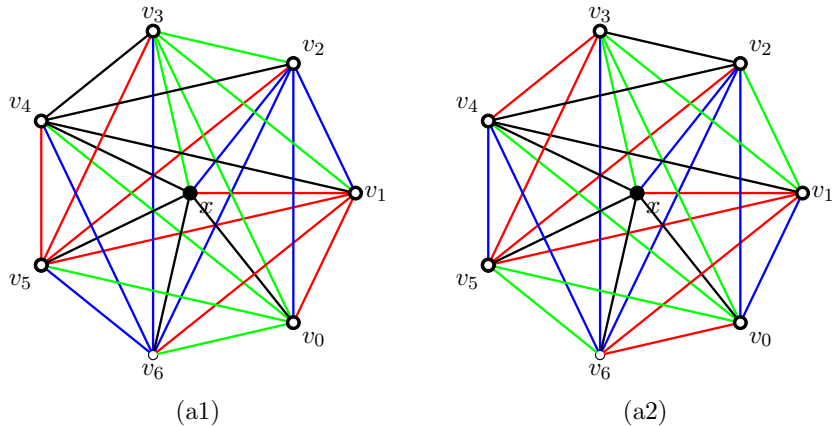

\vspace{1mm} {\em Case (b) $v_0v_{n+1} \in E(T_1)\/$  }

\vspace{1mm} By the argument in (i), $v_2v_{n} \in E(T_1)\/$.

\vspace{1mm} (b1) If $v_{2n-2}v_{n+1} \in E(T_1)\/$, then $v_{n-1}v_2 \in E(T_1)\/$. Moreover $v_{i-2}v_{n+i}, v_{n-i}v_{i+1} \in E(T_i)\/$, $i=1, 2, \ldots, n-1\/$.

\vspace{1mm} (b2) If $v_0v_{n+2} \in E(T_1)\/$, then $v_{n}v_3 \in E(T_1)\/$. Moreover $v_{i-1}v_{n+i+1}, v_{n+i-1}v_{i+2} \in E(T_i)\/$, $i=1, 2, \ldots, n-1\/$.

\vspace{1mm} Continue with this argument, we see that each of the plane trees $T_1, T_2, \ldots, T_{n-1}\/$ are $P_4\/$-symmetric of Type-2 and that they are pairwise isomorphic. Here the $4\/$-path in $T_i\/$ is $xv_iv_{n+i}v_{i+1}\/$ and  $\varphi(v_{i+n}) = v_{i+1}\/$. It is easy to see that the $T_i\/$'s in (b1) are caterpillars while the $T_i\/$'s in (b2) are such that $T_i -x\/$ is a caterpillar for each $i =1, 2, \ldots, n-1\/$.

\vspace{1mm} Hence, the $T_i\/$'s in two cases (a) and (b) are $w\/$-caterpillars.

\vspace{1mm} It is easy to see that for each non-boundary edge $v_iv_j \in E(T_n)\/$ exactly one of $v_iv_{j-1}\/$ and $v_jv_{i+1}\/$ is an edge of $T_n\/$, it follows that $T_n\/$ has only one boundary edge that lie on the segment from $v_1\/$ to $v_{n-1}\/$ (in anti-clockwise direction) since a non-boundary edge $v_1v_{n-1} \in E(T_n)\/$. Thus $T_n\/$ is a plane spanning caterpillar.

\begin{figure}[htb]
\centering
\resizebox{12cm}{!}{
\begin{minipage}{.45\textwidth}
\begin{tikzpicture}

       \coordinate (x) at (0,0);\node at (x) {\textbullet};
       \coordinate (v0) at (1.56,-1.96);\filldraw[black] (v0) circle(2pt);\node[below right] at (v0) {${v_0}$};
       \coordinate (v1) at (2.5,0);\filldraw[black] (v1) circle(2pt);\node[right] at (v1) {${v_1}$};
       \coordinate (v2) at (1.56,1.96);\filldraw[black] (v2) circle(2pt);\node[above right] at (v2) {${v_2}$};
       \coordinate (v3) at (-0.56,2.45);\filldraw[black] (v3) circle(2pt);\node[above] at (v3) {${v_3}$};
       \coordinate (v4) at (-2.25,1.09);\filldraw[black] (v4) circle(2pt);\node[above left] at (v4) {${v_4}$};
       \coordinate (v5) at (-2.25,-1.09);\filldraw[black] (v5) circle(2pt);\node[below left] at (v5) {${v_5}$};
       \coordinate (v6) at (-0.56,-2.45);\filldraw[black] (v6) circle(2pt);\node[below] at (v6) {${v_6}$};

       \draw [line width=1,red](v1) -- (x);\draw [line width=1,red](v1) -- (v5);\draw [line width=1,red](v5) -- (v2);
       \draw [line width=1,blue](v2) -- (x);\draw [line width=1,blue](v2) -- (v6);\draw [line width=1,blue](v6) -- (v3);
       \draw [line width=1,green](v3) -- (x);\draw [line width=1,green](v3) -- (v0);\draw [line width=1,green](v0) -- (v4);
       \draw [line width=1,black](v4) -- (x);\draw [line width=1,black](v5) -- (x);\draw [line width=1,black](v6) -- (x);
       \draw [line width=1,black](v0) -- (x);\draw [line width=1,black](v1) -- (v4);

      \draw [line width=1,red](v0) -- (v5);\draw [line width=1,red](v2) -- (v4);
      \draw [line width=1,blue](v1) -- (v6);\draw [line width=1,blue](v3) -- (v5);
      \draw [line width=1,green](v0) -- (v2);\draw [line width=1,green](v4) -- (v6);
      \draw [line width=1,black](v1) -- (v3);

      \draw [line width=1,red](v2) -- (v3);\draw [line width=1,red](v5) -- (v6);
      \draw [line width=1,blue](v3) -- (v4);\draw [line width=1,blue](v6) -- (v0);
      \draw [line width=1,green](v4) -- (v5);\draw [line width=1,green](v0) -- (v1);
      \draw [line width=1,black](v1) -- (v2);

       \filldraw[black] (v0) circle(3pt);\filldraw[black] (v1) circle(3pt);\filldraw[black] (v2) circle(3pt);
       \filldraw[black] (v3) circle(3pt);\filldraw[black] (v4) circle(3pt);\filldraw[black] (v5) circle(3pt);\filldraw[black] (v5) circle(3pt);

       \filldraw[white] (v0) circle(1.5pt);\filldraw[white] (v1) circle(1.5pt);\filldraw[white] (v2) circle(1.5pt);
       \filldraw[white] (v3) circle(1.5pt);\filldraw[white] (v4) circle(1.5pt);\filldraw[white] (v5) circle(1.5pt);\filldraw[white] (v6) circle(1.5pt);
       \filldraw[black] (x) circle(3pt);\node[below right] at (x) {${x}$};

\end{tikzpicture}
\centering

(b1)

\end{minipage}
\begin{minipage}{.45\textwidth}
\begin{tikzpicture}

       \coordinate (x) at (0,0);\node at (x) {\textbullet};
       \coordinate (v0) at (1.56,-1.96);\filldraw[black] (v0) circle(2pt);\node[below right] at (v0) {${v_0}$};
       \coordinate (v1) at (2.5,0);\filldraw[black] (v1) circle(2pt);\node[right] at (v1) {${v_1}$};
       \coordinate (v2) at (1.56,1.96);\filldraw[black] (v2) circle(2pt);\node[above right] at (v2) {${v_2}$};
       \coordinate (v3) at (-0.56,2.45);\filldraw[black] (v3) circle(2pt);\node[above] at (v3) {${v_3}$};
       \coordinate (v4) at (-2.25,1.09);\filldraw[black] (v4) circle(2pt);\node[above left] at (v4) {${v_4}$};
       \coordinate (v5) at (-2.25,-1.09);\filldraw[black] (v5) circle(2pt);\node[below left] at (v5) {${v_5}$};
       \coordinate (v6) at (-0.56,-2.45);\filldraw[black] (v6) circle(2pt);\node[below] at (v6) {${v_6}$};

       \draw [line width=1,red](v1) -- (x);\draw [line width=1,red](v1) -- (v5);\draw [line width=1,red](v5) -- (v2);
       \draw [line width=1,blue](v2) -- (x);\draw [line width=1,blue](v2) -- (v6);\draw [line width=1,blue](v6) -- (v3);
       \draw [line width=1,green](v3) -- (x);\draw [line width=1,green](v3) -- (v0);\draw [line width=1,green](v0) -- (v4);
       \draw [line width=1,black](v4) -- (x);\draw [line width=1,black](v5) -- (x);\draw [line width=1,black](v6) -- (x);
       \draw [line width=1,black](v0) -- (x);\draw [line width=1,black](v1) -- (v4);

        \draw [line width=1,red](v0) -- (v5);\draw [line width=1,red](v2) -- (v4);
      \draw [line width=1,blue](v1) -- (v6);\draw [line width=1,blue](v3) -- (v5);
      \draw [line width=1,green](v0) -- (v2);\draw [line width=1,green](v4) -- (v6);
      \draw [line width=1,black](v1) -- (v3);

      \draw [line width=1,red](v3) -- (v4);\draw [line width=1,red](v0) -- (v6);
      \draw [line width=1,blue](v0) -- (v1);\draw [line width=1,blue](v4) -- (v5);
       \draw [line width=1,green](v1) -- (v2);\draw [line width=1,green](v5) -- (v6);
        \draw [line width=1,black](v2) -- (v3);

       \filldraw[black] (v0) circle(3pt);\filldraw[black] (v1) circle(3pt);\filldraw[black] (v2) circle(3pt);
       \filldraw[black] (v3) circle(3pt);\filldraw[black] (v4) circle(3pt);\filldraw[black] (v5) circle(3pt);\filldraw[black] (v5) circle(3pt);

       \filldraw[white] (v0) circle(1.5pt);\filldraw[white] (v1) circle(1.5pt);\filldraw[white] (v2) circle(1.5pt);
       \filldraw[white] (v3) circle(1.5pt);\filldraw[white] (v4) circle(1.5pt);\filldraw[white] (v5) circle(1.5pt);\filldraw[white] (v6) circle(1.5pt);
       \filldraw[black] (x) circle(3pt);\node[below right] at (x) {${x}$};

\end{tikzpicture}
\centering

(b2)

\end{minipage}
}

\vspace{5mm}

\caption{\small { An illustration of Case (b). } \label{fig-4-th1}}
\end{figure}
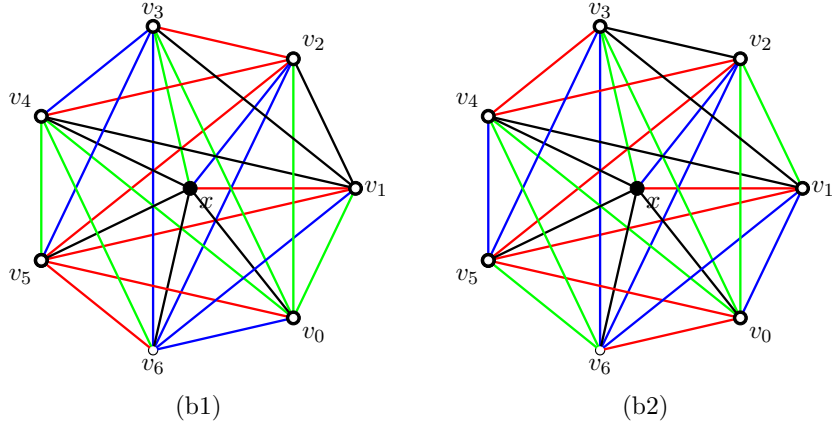

\vspace{1mm} Conversely, let $T\/$ be a $w\/$-caterpillar with $u_1u_2u_3u_4\/$ as the $4\/$-path.

\vspace{1mm} First place $u_1u_2u_3u_4\/$ on $P\/$ such that  $u_1 =x, u_2 =v_1, u_3 =v_{n+1}, u_4 = v_2\/$.

\vspace{1mm} Let $A, B\/$ be the components of $T- \{u_1u_2u_3u_4\}\/$. Since $T\/$ is a $w\/$-caterpillar, we see that $A\/$ and $B\/$ are isomorphic  caterpillars. We distinguish two cases.

\vspace{1mm} If $T\/$ is Type-1, we consider the segments $S_1 :v_{n+2}, v_{n+3}, \ldots, v_1\/$ and $S_2: v_3, v_4, \ldots, v_{n+1}\/$ (which are convex polygons). Place   $A\/$ and $B\/$ symmetrically (with respect to $u_2\/$ and $u_3\/$) on $S_1\/$ and $S_2\/$  respectively such that $S_1, S_2\/$ each has exactly $2\/$ boundary edges. Note that any caterpillar can be placed on a convex polygon with only two boundary edges. We omit the proof since it can be proved easily by induction.

\vspace{1mm}  If $T\/$ is Type-2, we consider the segments $S_1 :v_{n+1}, v_{n+2}, \ldots, v_0\/$ and $S_2: v_2, v_3, \ldots, v_{n}\/$ and place   $A\/$ and $B\/$ symmetrically (with respect to $u_3\/$ and $u_4\/$) on $S_1\/$ and $S_2\/$  respectively such that $S_1, S_2\/$ each has exactly $2\/$ boundary edges.

\vspace{1mm} Now let $T_1=T\/$. For each $i = 2, \ldots, n-1\/$, let $T_i \/$ be obtained by rotating the edges of $T_{i-1}\/$ with respect to the center vertex $x\/$ so that each non-center vertex of $T_{i-1}\/$ is rotated  anticlockwise exactly once on $C\/$.

\vspace{1mm} Clearly $xv_i \in E(T_n)\/$ for each $i=n, n+1, \ldots, 2n-2, 0\/$.

\vspace{1mm} Note that the edge $v_1v_n\/$ crosses an edge in the $4\/$-path $xv_iv_{i+n}v_{i+1}\/$ for each $i =1, \ldots, n-1\/$. This implies that  $v_1v_n \in E(T_n)\/$. It is clear that $T_n\/$ has only one boundary edge.

\vspace{1mm} There are two cases to consider.

\vspace{1mm} {\em Case (1):} $T_1\/$ is of Type-1.

\vspace{1mm} Then $v_1v_{n+2}, v_3v_{n+1} \in E(T_1)\/$ and  $v_1v_{n-1} \in E(T_{n-1})\/$ (and hence  $v_1v_{n-1}\notin E(T_n)\/$).


 \vspace{1mm} Since $v_2v_n\/$  crosses an edge of the $4\/$-path $xv_iv_{i+n}v_{i+1}\/$ for each $i =2, \ldots, n-1\/$, and $v_2v_n\/$ crosses $v_3v_{n+1}\/$, it follows that  $v_2v_n \in E(T_n)\/$.

\vspace{1mm} {\em Case (2):} $T_1\/$ is of Type-2.

Then $v_2v_{n}, v_0v_{n+1} \in E(T_1)\/$ (and hence  $v_2v_n\notin E(T_n)\/$).

 \vspace{1mm} Since $v_1v_{n-1}\/$  crosses an edge of the $4\/$-path $xv_iv_{i+n}v_{i+1}\/$ for each $i =1, 2, \ldots, n-2\/$, and $v_1v_{n-1}\/$ crosses $v_{n-2}v_0\/$, it follows that  $v_1v_{n-1} \in E(T_n)\/$.

\vspace{1mm}
 In either of the above cases, note that if $v_iv_j \in E(T_1)\/$  where $2 \leq i < j \leq n+1\/$,   then $v_{i+n-1}v_{j+n-1} \in E(T_1)\/$ and that  $v_{i-1}v_{j-1} \in E(T_n)\/$.

\vspace{1mm} This completes the proof. \qed

\begin{note}\label{rem1}
Note that in Theorem \ref{the1}, (i) if $T_1, T_2, \ldots, T_{n-1}\/$ are isomorphic double star, then $T_n\/$ is also a double star and $T_n \/$ is isomorphic to $T_i\/$ for all $i=1, 2, \ldots, n-1\/$. (ii) If $T_1, T_2, \ldots, T_{n-1}\/$ are isomorphic $w\/$-caterpillar but not double star, then $T_n \/$ is not isomorphic to $T_i\/$ for any  $i \in \{1, 2, \ldots, n-1\}\/$.
\end{note}

\section{Sufficient conditions}\label{sec-sc}

Given a set $P\/$ of $2m\/$  points in the plane in general position i.e., no three points are collinear. A {\em halving line\/} of $P\/$ is a line passing through two points in $P\/$, and cutting the remaining set of $2m-2\/$ points in half i.e., leaving $m-1\/$ points of $P\/$ on each side.

\vspace{1mm}
Given an arbitrary line passing through any point of $P\/$, we can always turn  it around by at most $180\/$ degrees where  it hits some other point of $P\/$  which makes it  a halving line. Thus the number of the halving lines is at least $m\/$.

\vspace{1mm}
In \cite{ps:refer} the authors proved that, a set of $2m\/$  points in general position in the plane admits a perfect crossing-matching (i.e., $m\/$ pairwise crossing edges) if and only if it has precisely $m\/$ halving lines.

\vspace{1mm}
Observe that for a set of $2m\/$ points in convex position in the plane, there are exactly $m\/$ halving lines. Hence, each of the $2m\/$ points is incident with exactly one halving line.

\begin{support}\label{lem6}
Let $P\/$ be a set of points in general position, and let $\ell\/$ be a halving line passing through $u,v \in P\/$. Let $P_1\/$ and $P_2\/$ be two the half-planes defined by $\ell\/$ and let $w_i \in P_i\/$, $i=1, 2\/$. Let $T\/$ be the geometric graph with vertex set $P\/$ and edge set $\{w_1u, uv,  vw_2\}\cup \{w_1x:x\in P_1\}\cup \{vy:y\in P_2\}\/$.  Then $T\/$ is a plane spanning tree.
\end{support}

\vspace{1mm}  \noindent
{\bf Proof:} The set of edges incident to $w_i\in P_i \/$ for $i\in\{1,2\} \/$ form a star. Regardless of the point set, a geometric star is always plane.  Thus no two edges incident to  $w_i\/$ for $i\in\{1,2\} \/$ cross each other.

\vspace{1mm} Since $P_1\cap P_2=  \emptyset\/$,  no edge incident to $w_i\/$  crosses $uv\/$ for each  $i= 1, 2\/$. Furthermore, no edge incident to $w_1\/$ crosses an edge incident to $w_2\/$ since such edges are separated by $\ell $. \qed

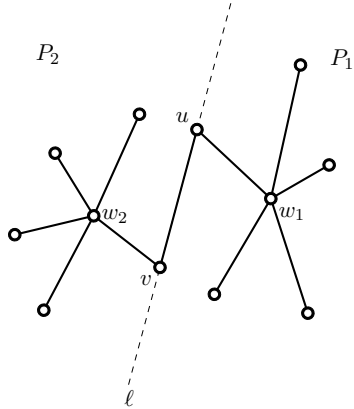
\begin{figure}[htb]
\centering
\resizebox{5cm}{!}{
\begin{tikzpicture}[rotate=75,.style={draw}]

       \coordinate (L1) at (-3.4,0);\node[below] at (L1) {${\ell}$};

       \coordinate (u) at (1,0);\node at (u) {\textbullet};\node[above left] at (u) {${u}$};
       \coordinate (v2) at (1,1);\node at (v2) {\textbullet};
       \coordinate (v3) at (-1.5,2.5);\node at (v3) {\textbullet};
       \coordinate (v4) at (0,2.2);\node at (v4) {\textbullet};
       \coordinate (w2) at (-0.85,1.3);\node at (w2) {\textbullet};\node[right] at (w2) {${w_2}$};
       \coordinate (v6) at (-2.6,1.7);\node at (v6) {\textbullet};
       \coordinate (v) at (-1.4,0);\node at (v) {\textbullet};\node[below left] at (v) {${v}$};
       \coordinate (v8) at (-1.6,-1);\node at (v8) {\textbullet};
       \coordinate (v9) at (-1.5,-2.6);\node at (v9) {\textbullet};
       \coordinate (w1) at (0.2,-1.5);\node at (w1) {\textbullet};\node[below right] at (w1) {${w_{1}}$};
       \coordinate (v11) at (1,-2.3);\node at (v11) {\textbullet};
       \coordinate (v0) at (2.5,-1.4);\node at (v0) {\textbullet};

        \coordinate (P1) at (2.4,-1.8);\node at (P1) {};\node[above right] at (P1) {${P_{1}}$};
         \coordinate (P2) at (1.2,3);\node at (P2) {};\node[above right] at (P2) {${P_{2}}$};

       \draw [line width=1,black](u) -- (v);\draw [line width=1,black](w1) -- (u);\draw [line width=1,black](w1) -- (v0);
       \draw [line width=1,black](w1) -- (v11);\draw [line width=1,black](w1) -- (v9);\draw [line width=1,black](w1) -- (v8);
       \draw [line width=1,black](w2) -- (v);\draw [line width=1,black](w2) -- (v2);\draw [line width=1,black](w2) -- (v3);
       \draw [line width=1,black](w2) -- (v4);\draw [line width=1,black](w2) -- (v6);
       \draw [dashed,black](3.2,0) --  (-3.5,0);
       \filldraw[black] (v0) circle(3pt); \filldraw[black] (v2) circle(3pt); \filldraw[black] (v4) circle(3pt); \filldraw[black] (v6) circle(3pt);
       \filldraw[black] (v8) circle(3pt); \filldraw[black] (w1) circle(3pt);
       \filldraw[black] (u) circle(3pt);\filldraw[black] (v3) circle(3pt);\filldraw[black] (w2) circle(3pt);
       \filldraw[black] (v) circle(3pt);\filldraw[black] (v9) circle(3pt);\filldraw[black] (v11) circle(3pt);
       \filldraw[white] (v0) circle(1.5pt); \filldraw[white] (v2) circle(1.5pt); \filldraw[white] (v4) circle(1.5pt);
       \filldraw[white] (v6) circle(1.5pt);
       \filldraw[white] (v8) circle(1.5pt); \filldraw[white] (w1) circle(1.5pt);
       \filldraw[white] (u) circle(1.5pt);\filldraw[white] (v3) circle(1.5pt);\filldraw[white] (w2) circle(1.5pt);
       \filldraw[white] (v) circle(1.5pt);\filldraw[white] (v9) circle(1.5pt);\filldraw[white] (v11) circle(1.5pt);
\end{tikzpicture}
}
\caption{An illustration of Lemma \ref{lem6} } \label{sufficint}

\end{figure}

\begin{support}\label{lem7}
Let $P\/$ be a set of points in general position. Let $\ell_1\/$ be a halving line passing through $u_1, u_2 \in P\/$ with $P_1\/$, $P_2\/$ the two half-planes defined by $\ell_1\/$. Suppose $v_1\in P_1\/$ and $v_2\in P_2\/$. Let $\ell_2\/$ be a halving line passing through $v_1\/$ and $v_2\/$ with $S_1\/$, $S_2\/$ the two half-planes defined by $\ell_2\/$ such that $u_1\in S_1\/$ and $u_2\in S_2\/$. Suppose $T_1\/$ and $T_2\/$ are the plane geometric graphs with vertex set $P\/$ and edge sets $\{u_1u_2,u_1v_2,u_2v_1\}\cup \{v_1x:x\in P_{1}\}  \cup \{v_2x:x\in P_{2}\} \/$ and   $\{v_1v_2,v_1u_1,v_2u_2\}\cup \{u_1y:y\in S_{1}\} \cup \{u_2y:y\in S_{2}\}  \/$ respectively. Then $E(T_1)\cap E(T_2)= \emptyset \/$.
\end{support}

\vspace{1mm}  \noindent
{\bf Proof:}  Assume on the contrary that there is an edge $e \in E(T_1)\cap E(T_2)\/$.

\vspace{1mm} Note that all edges of $T_1\/$ (but $u_1u_2\/$) are incident to $v_1\/$ or $v_2\/$,  and all edges of $T_2\/$ (but $v_1v_2\/$) are incident to $u_1\/$ or $u_2\/$. Hence $e\/$ must be one of the edges in the complete graph $K_4\/$ induced the vertices $\{u_1, u_2, v_1, v_2\}\/$.

\vspace{1mm} However  this is a contradiction because $E(T_1) \cap E(K_4) = \{ v_1u_2, u_2u_1, u_1v_2\}\/$ and  $E(T_2) \cap E(K_4) = \{ u_1v_1, v_1v_2, v_2u_2\}\/$.   Thus $E(T_1)\cap E(T_2)= \emptyset\/$. \qed

\vspace{1mm}
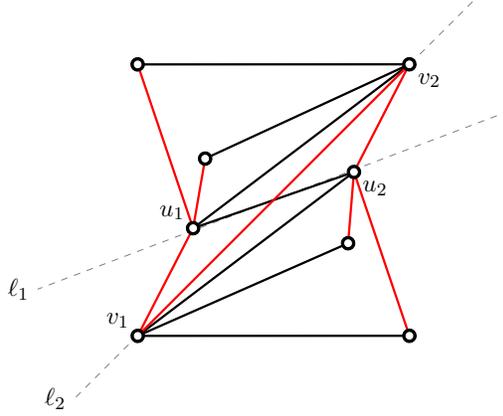
\begin{figure}[htb]
\centering
\resizebox{7cm}{!}{
\begin{tikzpicture}[rotate=75,.style={draw}]
\coordinate (center) at (0,0);
  \def\radius{3.7cm}

       \draw [dashed,gray](2.3,-3.2) -- (-2.3,3.2);
       \draw [dashed,gray](3.8,-2.2) -- (-3.8,2.2);
       \coordinate (L1) at (-2.3,3.2);\node[left] at (L1) {${\ell_1}$};
       \coordinate (L2) at (-3.8,2.2);\node[left] at (L2) {${\ell_2}$};
        \coordinate (u1) at (-0.75,1.1);\node at (u1) {\textbullet};\node[above left] at (u1) {${u_1}$};
        \coordinate (u2) at (0.75,-1.1);\node at (u2) {\textbullet};\node[below right] at (u2) {${u_{2}}$};
        \coordinate (v1) at (-2.6,1.5);\node at (v1) {\textbullet};\node[above left] at (v1) {${v_1}$};
        \coordinate (v2) at (2.6,-1.5);\node at (v2) {\textbullet};\node[below right] at (v2) {${v_{2}}$};
        \coordinate (w1) at (-1.5,-2.6);
        \coordinate (w2) at (1.5,2.6);
        \coordinate (z1) at (-0.35,-1.3);
        \coordinate (z2) at (0.35,1.2);
      \draw [line width=1,black](u1) -- (u2);\draw [line width=1,black](u1) -- (v2);\draw [line width=1,black](u2) -- (v1);
      \draw [line width=1,black](v1) -- (z1);\draw [line width=1,black](v1) -- (w1);\draw [line width=1,black](v2) -- (z2);
      \draw [line width=1,black](v2) -- (w2);
      \draw [line width=1,red](v1) -- (v2);\draw [line width=1,red](v1) -- (u1); \draw [line width=1,red](v2) -- (u2);
      \draw [line width=1,red](u1) -- (w2); \draw [line width=1,red](u1) -- (z2); \draw [line width=1,red](u2) -- (w1);
      \draw [line width=1,red](u2) -- (z1);
       \filldraw[black] (w2) circle(3pt);
       \filldraw[black] (v2) circle(3pt);
       \filldraw[black] (v1) circle(3pt);
       \filldraw[black] (z1) circle(3pt);
       \filldraw[black] (u2) circle(3pt);
       \filldraw[black] (w1) circle(3pt);
       \filldraw[black] (z2) circle(3pt);
       \filldraw[black] (u1) circle(3pt);
       \filldraw[white] (w2) circle(1.5pt);
       \filldraw[white] (v2) circle(1.5pt);
       \filldraw[white] (v1) circle(1.5pt);
       \filldraw[white] (z1) circle(1.5pt);
       \filldraw[white] (u2) circle(1.5pt);
       \filldraw[white] (w1) circle(1.5pt);
       \filldraw[white] (z2) circle(1.5pt);
       \filldraw[white] (u1) circle(1.5pt);
\end{tikzpicture}
}
\caption{An illustration of  Lemma \ref{lem7}} \label{sufficint}
\end{figure}

\begin{define}\label{def3}
Recall that  $S(m, n)\/$ is the double star with two non-pendant vertices  $v, w\/$. Let $S_k(m,n)\/$ be the tree obtained from $S(m,n)\/$ by inserting $k\/$ vertices of degree $2\/$ into the edge $vw\/$. Clearly $S_0(m,n)\/$ is the double star $S(m,n)\/$.
\end{define}

\vspace{1mm}
A set $P\/$ of $2m\/$ points is said to have a {\em $h\/$-labeling\/} if it has precisely $m\/$ halving lines $h_0, h_1, \ldots, h_{m-1}\/$ in anticlockwise direction.  In this case, we label the two points of $P\/$ that are incident to $h_i\/$ with $v_i, v_{m+i}\/$ so that $v_0, v_1, \ldots, v_{2m-1}\/$ are in anticlockwise direction.

\begin{result}\label{the2}
Let $P\/$ be a set of $2n\/$ points in general position, $n \geq 2\/$. Suppose each point of $P\/$ is incident with exactly one halving line. Then the edges of the complete geometric graph on $P\/$ can be partitioned into $n\/$ plane spanning trees $S_{2t}(m,m)\/$ where $t\leq 1\/$ and $m = n-t\/$.
\end{result}

\vspace{1mm}  \noindent
{\bf Proof:} The case $t=0\/$ has been treated in \cite{phrw:refer}. Thus we assume that $t =1\/$.

\vspace{1mm}
Since each point of $P\/$ is incident with exactly one halving line, it follows that $P\/$ has a $h\/$-labeling $h_i\/$,  $i=0, \ldots, n-1\/$ where $h_i\/$ passes through $v_i, v_{i+n}\/$.

\vspace{1mm} Note that each halving line $h_i\/$ splits $P\/$ into two  point sets $P_{i,1} = \{ v_{i+1}, v_{i+2}, \ldots, v_{i+n-1}\}\/$ and $P_{i,2}  = \{v_{i+n+1}, v_{i+n+2}, \ldots, v_{i+2n-1}\}\/$.

\vspace{1mm} Let $T_i=\{v_iv_{i+n},v_{i}v_{i+n+1},v_{i+n}v_{i+1}\}\cup \{v_{i+1}x:x\in P_{i,1} \}\cup \{v_{i+n+1}y:y\in P_{i,2} \}\/$ for $i=0,\dots, n-2\/$.
Then $T_i\/$ is the tree $S_2(n-1,n-1)\/$  on $P\/$. Here the operations on the subscripts are reduced modulo $2n\/$.

\vspace{1mm} By Lemma \ref{lem6} each $T_i\/$ is plane tree. Moreover, by Lemma \ref{lem7},  the $T_i\/$'s are pairwise edge-disjoint. \qed

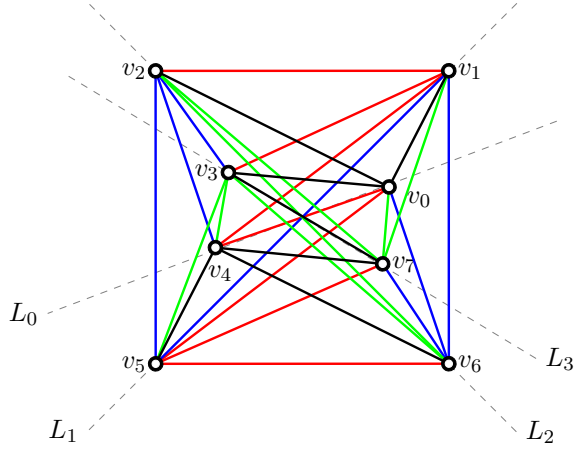
\begin{figure}[htb]
\centering
\resizebox{8cm}{!}{
\begin{tikzpicture}[rotate=75,.style={draw}]
\coordinate (center) at (0,0);
  \def\radius{3.7cm}

       \draw [dashed,gray](2.3,-3.2) -- (-2.3,3.2);
       \draw [dashed,gray](3.8,-2.2) -- (-3.8,2.2);
       \draw [dashed,gray](-2.2,-3.8) -- (2.2,3.8);
       \draw [dashed,gray](-1.1,-3.8) -- (1.1,3.8);
       \coordinate (L0) at (-2.3,3.2);\node[left] at (L0) {${L_0}$};
       \coordinate (L1) at (-3.8,2.2);\node[left] at (L1) {${L_1}$};
       \coordinate (L2) at (-2.2,-3.8);\node[right] at (L2) {${L_2}$};
       \coordinate (L3) at (-1.1,-3.8);\node[right] at (L3) {${L_3}$};
       \coordinate (v2) at (1.5,2.6);\node at (v2) {\textbullet};\node[left] at (v2) {${v_2}$};
       \coordinate (v3) at (0.35,1.2);\node at (v3) {\textbullet};\node[left] at (v3) {${v_3}$};
       \coordinate (v4) at (-0.75,1.1);\node at (v4) {\textbullet};\node[below] at (-0.85,1) {${v_4}$};
       \coordinate (v5) at (-2.6,1.5);\node at (v5) {\textbullet};\node[left] at (v5) {${v_5}$};
       \coordinate (v7) at (-0.35,-1.3);\node at (v7) {\textbullet};\node[right] at (v7) {${v_7}$};
       \coordinate (v6) at (-1.5,-2.6);\node at (v6) {\textbullet};\node[right] at (v6) {${v_6}$};
       \coordinate (v0) at (0.75,-1.1);\node at (v0) {\textbullet};\node[below right] at (0.85,-1.2) {${v_{0}}$};
       \coordinate (v1) at (2.6,-1.5);\node at (v1) {\textbullet};\node[right] at (v1) {${v_{1}}$};

      \draw [line width=1,red](v4) -- (v0);\draw [line width=1,red](v4) -- (v1);\draw [line width=1,red](v0) -- (v5);
      \draw [line width=1,red](v5) -- (v7);\draw [line width=1,red](v5) -- (v6);\draw [line width=1,red](v1) -- (v3);
      \draw [line width=1,red](v1) -- (v2);

    \draw [line width=1,blue](v5) -- (v1);\draw [line width=1,blue](v5) -- (v2); \draw [line width=1,blue](v1) -- (v6);
     \draw [line width=1,blue](v2) -- (v4);\draw [line width=1,blue](v2) -- (v3);
     \draw [line width=1,blue](v6) -- (v0);\draw [line width=1,blue](v6) -- (v7);

      \draw [line width=1,green](v2) -- (v6); \draw [line width=1,green](v6) -- (v3); \draw [line width=1,green](v2) -- (v7);\draw [line width=1,green](v7) -- (v1); \draw [line width=1,green](v7) -- (v0);\draw [line width=1,green](v3) -- (v5); \draw [line width=1,green](v3) -- (v4);

       \draw [line width=1,black](v3) -- (v7); \draw [line width=1,black](v3) -- (v0); \draw [line width=1,black](v7) -- (v4); \draw [line width=1,black](v4) -- (v5); \draw [line width=1,black](v4) -- (v6);\draw [line width=1,black](v0) -- (v1); \draw [line width=1,black](v0) -- (v2);
       \filldraw[black] (v2) circle(3pt);
       \filldraw[black] (v1) circle(3pt);
       \filldraw[black] (v5) circle(3pt);
       \filldraw[black] (v7) circle(3pt);
       \filldraw[black] (v0) circle(3pt);
       \filldraw[black] (v6) circle(3pt);
       \filldraw[black] (v3) circle(3pt);
       \filldraw[black] (v4) circle(3pt);
       \filldraw[white] (v2) circle(1.5pt);
       \filldraw[white] (v1) circle(1.5pt);
       \filldraw[white] (v5) circle(1.5pt);
       \filldraw[white] (v7) circle(1.5pt);
       \filldraw[white] (v0) circle(1.5pt);
       \filldraw[white] (v6) circle(1.5pt);
       \filldraw[white] (v3) circle(1.5pt);
       \filldraw[white] (v4) circle(1.5pt);
\end{tikzpicture}
}
\caption{An illustration of  Theorem \ref{the2}}
\label{fig4the2}
\end{figure}

Let $P\/$ be a set of $2m\/$ points in general position in the plane, $ m\geq 2\/$. We say that $L\/$ is a $k\/$-halving line of $P\/$ if $L\/$ is a line passing through two points of $P\/$ and splits the remaining $2m-2\/$ points into two sets one of size $m-1-k\/$ and another of size $m-1+k\/$ where $k=0, 1, \ldots, m-1\/$. A $k\/$-halving line that passes through  to the two points $u\/$ and $v\/$ will be denoted by $h^{k}(u,v)\/$. Clearly $h^0(u,v)\/$ is a halving line of $P\/$.

\begin{result}\label{the3}
Let $P\/$ be a set of $2n\/$ points in general position, $n \geq 2\/$. Let $r\/$ be a natural number such that $1 \leq r \leq n-1\/$. Suppose   each point of $P\/$ is incident to exactly $b\/$ $k\/$-halving lines for every $ k=1, 2, \ldots, r\/$. Here $b = 1\/$ if $k=0\/$ and $b =2\/$ otherwise.   Then the edges of the complete geometric graph on $P\/$ can be partitioned into $n\/$ plane spanning symmetric caterpillars (that are pairwise graph-isomorphic).
\end{result}

\vspace{1mm}  \noindent
{\bf Proof:}  Since  each point is incident to exactly one halving line, it follows that $P\/$ has a $h\/$-labeling $h_i\/$,  $i=0, \ldots, n-1\/$  (where $h_i\/$ passes through $v_i, v_{i+n}\/$).

\vspace{1mm} Note that each halving line $h_i\/$ splits $P - \{v_i, v_{i+n}\}\/$ into two  point sets $P_{i,1} = \{ v_{i+1}, v_{i+2}, \ldots, v_{i+n-1}\}\/$ and $P_{i,2}  = \{v_{i+n+1}, v_{i+n+2}, \ldots, v_{i+2n-1}\}\/$.

\vspace{1mm}
Before going into the construction plane spanning trees on $P\/$, we first make the following observation.

\vspace{1mm} Let $L_1\/$ be the line passing through $v_i, v_{i+n-1}\/$. Suppose $L_1\/$ splits  $P-\{v_i, v_{i+n-1}\}\/$ into two point sets $P_1\/$ and $ P_2\/$ where $|P_1|<|P_2|\/$.

\vspace{1mm} We assert that $L_1\/$ is a $1\/$-halving line.

\vspace{1mm} To see this, we first note that  $P_1 \subseteq \{ v_{i+1}, v_{i+2}, \ldots, v_{i+n-2}\}\/$ since $h_i\/$ is a halving line with $P_{i,1} = \{ v_{i+1}, v_{i+2}, \ldots, v_{i+n-1}\}\/$.

\vspace{1mm}
If  $P_1 =\{ v_{i+1}, v_{i+2}, \ldots, v_{i+n-2}\}\/$, then the assertion follows. Hence we assume on the contrary that for some $1 \leq j \leq n-2\/$, $v_{i+j} \not \in P_1\/$.   Moreover $v_{i+j}\/$ can be chosen such that  the line $\ell_1\/$ passing through $v_i, v_{i+j}\/$ is a $1\/$-halving line.%

\vspace{1mm} Let $\ell_2\/$ be the line obtained from $\ell_1\/$ by fixing  the point $v_{i+j}\/$ and turn it around in anticlockwise direction until it becomes a $1\/$-halving line. Then either $\ell_2\/$ contains $v_{i+n}\/$ or else $\ell_2\/$ contains a vertex $v_{i+s} \/$ where $1 \leq s < j\/$.

\vspace{1mm} Let $\ell_3\/$ be the line obtained from $\ell_1\/$ by fixing  the point $v_{i+j}\/$ and turn it around in clockwise direction until it becomes a $1\/$-halving line. Then either $\ell_3\/$ contains $v_{i+n-1}\/$ or else $\ell_3\/$ contains a vertex $v_{i+t} \/$ where $j < t < n-1\/$.

\vspace{1mm} But this means that we have three $1\/$-halving lines passing through $v_{i+j}\/$,  a contradiction and the assertion follows.

\vspace{1mm} More generally, let $L_k\/$ be the line passing through $v_i, v_{i+n-k}\/$. Suppose $L_k\/$ splits  $P-\{v_i, v_{i+n-k}\}\/$ into two point sets $S_1\/$ and $ S_2\/$ where $|S_1|<|S_2|\/$.

\vspace{1mm} We assert that $L_k\/$ is a $k\/$-halving line if $k \leq r\/$

\vspace{1mm} We shall prove this assertion by induction on $k\/$. Since the assertion is true for $k=1\/$, we assume that $k \geq 2\/$.

\vspace{1mm} To see this, we first note that  $S_1 \subseteq \{ v_{i+1}, v_{i+2}, \ldots, v_{i+n-k-1}\}\/$ since  $L_{k-1}\/$ is the $(k-1)$-halving line which splits $P-\{v_i, v_{i+n-k+1}\}\/$ into two point sets   $\{v_{i+1}, v_{i+2}, \ldots, v_{i+n-k}\}\/$ and $ \{v_{i+n-k+2}, v_{i+n-k+3}, \ldots, v_{i+2n-1}\}\/$.

\vspace{1mm} If  $S_1 =\{ v_{i+1}, v_{i+2}, \ldots, v_{i+n-k-1}\}\/$, then the assertion follows.  Hence we assume on the contrary that for some $1 \leq j \leq n-k-1\/$, $v_{i+j} \not \in S_1\/$.   Moreover $v_{i+j}\/$ can be chosen such that  the line $\ell_1\/$ passing through $v_i, v_{i+j}\/$ is a $(k-1)\/$-halving line.

\vspace{1mm} By using an  argument similar to  the case $k=1\/$, we can show that there are three $(k-1)\/$-halving lines passing through $v_{i+j}\/$ giving a contradiction and the assertion follows.

\vspace{1mm}
Let $T_0\/$ be a tree on $P\/$ constructed in the following way.

\vspace{1mm}
(i) Let  $v_0v_n \in E(T_0)\/$.

\vspace{1mm}
(ii) Let $v_s,v_t \in P_{0,1}\cup \{v_0\} \/$ be such that the line passing through $v_s,v_t\/$ is the $r\/$-halving line $h^{r}(v_s,v_t)\/$ where $0 \leq s < t \leq n-1\/$. Then either choose $v_s x\/$ to be  an edge of  $T_0\/$ for every $x \in \{v_{s+1}, v_{s+2}, \ldots, v_{t}\} \/$,  or else choose $v_t y\/$ to be an edge of $T_0\/$ for every   $ y \in \{v_{s}, v_{s+1}, \ldots, v_{t-1} \}\/$.

\vspace{1mm}
(iii) With respect to the edge $v_sv_t\/$ of $T_0\/$, we either choose $v_tv_{s-1}\/$ (where $v_{s-1} \not \in P_{0,2}\/$)  or else $v_sv_{t+1}\/$  to be an edge of $T_0\/$. Iteratively, if $v_pv_q\/$ has been chosen to be an edge of $T_0\/$, where $0 \leq p \leq s\/$, $t \leq q \leq n-1\/$, then  we either choose $v_qv_{p-1}\/$ (where $v_{p-1} \not \in P_{0,2}\/$)    or else   $v_pv_{q+1}\/$  to be an edge of $T_0\/$.

\vspace{1mm}
(iv) Now to  every edge $v_iv_j \/$ that has been chosen to be in  $T_0\/$ where $i, j \in \{0, 1, \ldots, n-1\}\/$, choose $v_{i+n}v_{j+n}\/$ to be an edge of $T_0\/$.

\begin{figure}[htb]
\resizebox{15cm}{!}{
\begin{minipage}{.45\textwidth}
\begin{tikzpicture}[rotate=75,.style={draw}]
\coordinate (center) at (0,0);
  \def\radius{4cm}
       \coordinate (v1) at (1.9,0);\node at (v1) {\textbullet};\node[above] at (v1) {${v_1}$};
       \coordinate (v2) at (1.3,0.75);\node at (v2) {\textbullet};\node[above] at (v2) {${v_2}$};
       \coordinate (v3) at (1.5,2.6);\node at (v3) {\textbullet};\node[above left] at (v3) {${v_3}$};
       \coordinate (v4) at (0,1.9);\node at (v4) {\textbullet};\node[left] at (v4) {${v_4}$};
       \coordinate (v5) at (-0.75,1.3);\node at (v5) {\textbullet};\node[left] at (v5) {${v_5}$};
       \coordinate (v6) at (-2.6,1.5);\node at (v6) {\textbullet};\node[below left] at (v6) {${v_6}$};
       \coordinate (v7) at (-1.9,0);\node at (v7) {\textbullet};\node[below] at (v7) {${v_7}$};
       \coordinate (v8) at (-1.3,-0.75);\node at (v8) {\textbullet};\node[below] at (v8) {${v_8}$};
       \coordinate (v9) at (-1.5,-2.6);\node at (v9) {\textbullet};\node[below right] at (v9) {${v_9}$};
       \coordinate (v10) at (0,-1.9);\node at (v10) {\textbullet};\node[right] at (v10) {${v_{10}}$};
       \coordinate (v11) at (0.75,-1.3);\node at (v11) {\textbullet};\node[right] at (v11) {${v_{11}}$};
       \coordinate (v0) at (2.6,-1.5);\node at (v0) {\textbullet};\node[above right] at (v0) {${v_{0}}$};

       \draw [line width=1,red](v0) -- (v6);\draw [line width=1,red](v6) -- (v1);\draw [line width=1,red](v1) -- (v5);
       \draw [line width=1,red](v5) -- (v2);\draw [line width=1,red](v5) -- (v4);\draw [line width=1,red](v5) -- (v3);
       \draw [line width=1,red](v0) -- (v7);\draw [line width=1,red](v7) -- (v11);\draw [line width=1,red](v11) -- (v8);
       \draw [line width=1,red](v11) -- (v10);\draw [line width=1,red](v11) -- (v9);

       \filldraw[black] (v0) circle(3pt); \filldraw[black] (v2) circle(3pt); \filldraw[black] (v4) circle(3pt); \filldraw[black] (v6) circle(3pt);
       \filldraw[black] (v8) circle(3pt); \filldraw[black] (v10) circle(3pt);
       \filldraw[black] (v1) circle(3pt);\filldraw[black] (v3) circle(3pt);\filldraw[black] (v5) circle(3pt);\filldraw[black] (v5) circle(3pt);
       \filldraw[black] (v7) circle(3pt);\filldraw[black] (v9) circle(3pt);\filldraw[black] (v11) circle(3pt);

        \filldraw[white] (v0) circle(1.5pt); \filldraw[white] (v2) circle(1.5pt); \filldraw[white] (v4) circle(1.5pt); \filldraw[white] (v6) circle(1.5pt);
       \filldraw[white] (v8) circle(1.5pt); \filldraw[white] (v10) circle(1.5pt);
       \filldraw[white] (v1) circle(1.5pt);\filldraw[white] (v3) circle(1.5pt);\filldraw[white] (v5) circle(1.5pt);\filldraw[white] (v5) circle(1.5pt);
       \filldraw[white] (v7) circle(1.5pt);\filldraw[white] (v9) circle(1.5pt);\filldraw[white] (v11) circle(1.5pt);

\end{tikzpicture}
\centering

(a) \ $T_0$
\end{minipage}
\begin{minipage}{.45\textwidth}
\begin{tikzpicture}[rotate=75,.style={draw}]
\coordinate (center) at (0,0);
  \def\radius{4cm}
       \coordinate (v1) at (1.9,0);\node at (v1) {\textbullet};\node[above] at (v1) {${v_1}$};
       \coordinate (v2) at (1.3,0.75);\node at (v2) {\textbullet};\node[above] at (v2) {${v_2}$};
       \coordinate (v3) at (1.5,2.6);\node at (v3) {\textbullet};\node[above left] at (v3) {${v_3}$};
       \coordinate (v4) at (0,1.9);\node at (v4) {\textbullet};\node[above] at (v4) {${v_4}$};
       \coordinate (v5) at (-0.75,1.3);\node at (v5) {\textbullet};\node[above] at (v5) {${v_5}$};
       \coordinate (v6) at (-2.6,1.5);\node at (v6) {\textbullet};\node[below left] at (v6) {${v_6}$};
       \coordinate (v7) at (-1.9,0);\node at (v7) {\textbullet};\node[below] at (v7) {${v_7}$};
       \coordinate (v8) at (-1.3,-0.75);\node at (v8) {\textbullet};\node[below] at (v8) {${v_8}$};
       \coordinate (v9) at (-1.5,-2.6);\node at (v9) {\textbullet};\node[below right] at (v9) {${v_9}$};
       \coordinate (v10) at (0,-1.9);\node at (v10) {\textbullet};\node[below] at (v10) {${v_{10}}$};
       \coordinate (v11) at (0.75,-1.3);\node at (v11) {\textbullet};\node[below] at (v11) {${v_{11}}$};
       \coordinate (v0) at (2.6,-1.5);\node at (v0) {\textbullet};\node[above right] at (v0) {${v_{0}}$};

       \draw [line width=1,blue](v1) -- (v7);\draw [line width=1,blue](v7) -- (v2);\draw [line width=1,blue](v2) -- (v6);\draw [line width=1,blue](v6) -- (v3);\draw [line width=1,blue](v6) -- (v5);\draw [line width=1,blue](v6) -- (v4);\draw [line width=1,blue](v1) -- (v8);\draw [line width=1,blue](v8) -- (v0);\draw [line width=1,blue](v0) -- (v9);   \draw [line width=1,blue](v0) -- (v11);\draw [line width=1,blue](v0) -- (v10);

 \filldraw[black] (v0) circle(3pt); \filldraw[black] (v2) circle(3pt); \filldraw[black] (v4) circle(3pt); \filldraw[black] (v6) circle(3pt);
       \filldraw[black] (v8) circle(3pt); \filldraw[black] (v10) circle(3pt);
       \filldraw[black] (v1) circle(3pt);\filldraw[black] (v3) circle(3pt);\filldraw[black] (v5) circle(3pt);\filldraw[black] (v5) circle(3pt);
       \filldraw[black] (v7) circle(3pt);\filldraw[black] (v9) circle(3pt);\filldraw[black] (v11) circle(3pt);

        \filldraw[white] (v0) circle(1.5pt); \filldraw[white] (v2) circle(1.5pt); \filldraw[white] (v4) circle(1.5pt); \filldraw[white] (v6) circle(1.5pt);
       \filldraw[white] (v8) circle(1.5pt); \filldraw[white] (v10) circle(1.5pt);
       \filldraw[white] (v1) circle(1.5pt);\filldraw[white] (v3) circle(1.5pt);\filldraw[white] (v5) circle(1.5pt);\filldraw[white] (v5) circle(1.5pt);
       \filldraw[white] (v7) circle(1.5pt);\filldraw[white] (v9) circle(1.5pt);\filldraw[white] (v11) circle(1.5pt);

\end{tikzpicture}
\centering

(b) \ $T_1$
\end{minipage}
\begin{minipage}{.45\textwidth}
\begin{tikzpicture}[rotate=75,.style={draw}]
\coordinate (center) at (0,0);
  \def\radius{4cm}
       \coordinate (v1) at (1.9,0);\node at (v1) {\textbullet};\node[above] at (v1) {${v_1}$};
       \coordinate (v2) at (1.3,0.75);\node at (v2) {\textbullet};\node[above] at (v2) {${v_2}$};
       \coordinate (v3) at (1.5,2.6);\node at (v3) {\textbullet};\node[above left] at (v3) {${v_3}$};
       \coordinate (v4) at (0,1.9);\node at (v4) {\textbullet};\node[left] at (v4) {${v_4}$};
       \coordinate (v5) at (-0.75,1.3);\node at (v5) {\textbullet};\node[left] at (v5) {${v_5}$};
       \coordinate (v6) at (-2.6,1.5);\node at (v6) {\textbullet};\node[below left] at (v6) {${v_6}$};
       \coordinate (v7) at (-1.9,0);\node at (v7) {\textbullet};\node[below] at (v7) {${v_7}$};
       \coordinate (v8) at (-1.3,-0.75);\node at (v8) {\textbullet};\node[below] at (v8) {${v_8}$};
       \coordinate (v9) at (-1.5,-2.6);\node at (v9) {\textbullet};\node[below right] at (v9) {${v_9}$};
       \coordinate (v10) at (0,-1.9);\node at (v10) {\textbullet};\node[right] at (v10) {${v_{10}}$};
       \coordinate (v11) at (0.75,-1.3);\node at (v11) {\textbullet};\node[right] at (v11) {${v_{11}}$};
       \coordinate (v0) at (2.6,-1.5);\node at (v0) {\textbullet};\node[above right] at (v0) {${v_{0}}$};

       \draw [line width=1,black](v2) -- (v8);\draw [line width=1,black](v8) -- (v3);\draw [line width=1,black](v3) -- (v7);\draw [line width=1,black](v7) -- (v4);\draw [line width=1,black](v7) -- (v6);\draw [line width=1,black](v7) -- (v5);\draw [line width=1,black](v2) -- (v9);\draw [line width=1,black](v9) -- (v1);\draw [line width=1,black](v1) -- (v10);  \draw [line width=1,black](v1) -- (v0);\draw [line width=1,black](v1) -- (v11);

       \filldraw[black] (v0) circle(3pt); \filldraw[black] (v2) circle(3pt); \filldraw[black] (v4) circle(3pt); \filldraw[black] (v6) circle(3pt);
       \filldraw[black] (v8) circle(3pt); \filldraw[black] (v10) circle(3pt);
       \filldraw[black] (v1) circle(3pt);\filldraw[black] (v3) circle(3pt);\filldraw[black] (v5) circle(3pt);\filldraw[black] (v5) circle(3pt);
       \filldraw[black] (v7) circle(3pt);\filldraw[black] (v9) circle(3pt);\filldraw[black] (v11) circle(3pt);

        \filldraw[white] (v0) circle(1.5pt); \filldraw[white] (v2) circle(1.5pt); \filldraw[white] (v4) circle(1.5pt); \filldraw[white] (v6) circle(1.5pt);
       \filldraw[white] (v8) circle(1.5pt); \filldraw[white] (v10) circle(1.5pt);
       \filldraw[white] (v1) circle(1.5pt);\filldraw[white] (v3) circle(1.5pt);\filldraw[white] (v5) circle(1.5pt);\filldraw[white] (v5) circle(1.5pt);
       \filldraw[white] (v7) circle(1.5pt);\filldraw[white] (v9) circle(1.5pt);\filldraw[white] (v11) circle(1.5pt);

\end{tikzpicture}
\centering

(c) \ $T_2$
\end{minipage}
}
\caption{ An example of $T_0\/$  and its derived trees $T_1\/$ and $T_2\/$ in Theorem \ref{the3}. \label{t012}  }
\end{figure}
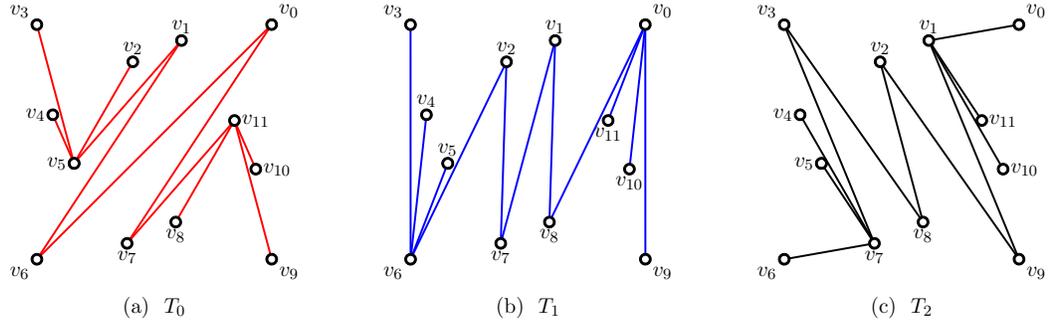

\vspace{1mm}
Clearly, $T_0\/$ is a spanning caterpillar on $P\/$. Moreover $T_0\/$ is a symmetric  because $T_0 - \{v_0v_n\}\/$ consists two isomorphic components $A, B\/$ where $A= P_{0,1} \cup \{v_0\}\/$,  $B = P_{0,2} \cup \{v_n\}\/$ and $v_0, v_n\/$ are similar in  $T_0 - \{v_0v_n\}\/$.  Figure \ref{t012} depicts an example of a spanning tree $T_0\/$ constructed in this manner.

\vspace{1mm} We now show that $T_0\/$ is a plane tree. Note that the subtree as constructed in (ii) is a  plane (because it is a geometric star). Also it follows from the observation (preceding the construction of $T_0\/$) that the subtree as constructed in (iii) is a plane tree. Moreover, no edge of $T_0\/$ from (ii)  crosses an edge of $T_0\/$ from (iii) since such edges are separated by the $r\/$-halving line $h^{r}(v_s,v_t)\/$. Furthermore, no edge of $T_0\/$ from (ii) and (iii) crosses an edge of $T_0\/$ from (iv) since such edges are separated by the halving line  $h_0\/$ (which  passes through $v_0, v_n\/$). We conclude that  $T_0\/$ is plane spanning symmetric caterpillar.

\vspace{1mm}
For each $q \in \{1, 2, \ldots, n-1\}\/$, let $T_q\/$ denote the tree on $P\/$ defined by $v_iv_j \in E(T_q)\/$ if and only if $v_{i-q}v_{j-q} \in E(T_0)\/$. Clearly $T_q\/$ is isomorphic to $T_0\/$ (see for example the trees $T_1\/$ and $T_2\/$ in Figure \ref{t012}).  It remains to show that  $E(T_0) \cap E(T_q) = \emptyset\/$.

\vspace{1mm}
Let $v_iv_j, v_lv_m \in E(T_0)\/$.  It follows from the above observation  and the definition of $T_0\/$ that if $i, j, l, m \in \{0,1, 2, \ldots, n-1\}\/$ where $i < j\/$ and $l < m\/$  (or $i, j, l, m \in \{n, n+1, n+2, \ldots, 2n-1\}\/$), then $|i-j| \neq |l-m|\/$. In the case that $|i-j| = |l-m|\/$, then $i = l+n\/$ and $j = m +n\/$ .

\vspace{1mm} Suppose $v_{i+q}v_{j+q} \in E(T_q) \cap E(T_0)\/$ (where $i, j \in \{0, 1, 2, \ldots, n-1\}\/$). By the definition of $T_0\/$, we have $v_iv_j, v_{i+q}v_{j+q} \in E(T_0)\/$. But then this implies that $q=n\/$, a contradiction. \qed

\vspace{1mm}
A set  $P\/$ of $2m\/$ points  is said to have a {\em $w\/$-labeling\/} if there is a vertex $w \in P\/$ incident to precisely $2m-1\/$ halving lines $h_0, h_1, \ldots, h_{2m-2}\/$ in anticlockwise direction. In this case, we label the points  of $P-\{w\}\/$ in anticlockwise direction $v_0, v_1, \ldots, v_{2m-2}\/$ so that $w\/$ and $v_i\/$ are incident to  $h_i\/$. Note that if $P\/$ is in a regular wheel configuration, then $P\/$ admits a $w\/$-labeling.

\begin{result}\label{the4}
Let $P\/$ be a set of $2n\/$ points in general position, $n \geq 2\/$. Let $r\/$ be a natural number such that $1 \leq r \leq n-1\/$. Suppose $w \in P\/$ is incident to $2n-1\/$ halving lines and each point of $P-\{w\}\/$ is incident to exactly two $k\/$-halving lines for every $ k=1, 2, \ldots, r\/$. Then the edges of the complete geometric graph on $P\/$ can be partitioned into $n\/$ plane spanning trees where $n-1\/$ of these trees are w-caterpillars
and are pairwise graph-isomorphic.
\end{result}

\vspace{1mm}  \noindent
{\bf Proof:}  The proof is similar to the proof of Theorem \ref{the3}. Hence we shall omit the details.

\vspace{1mm} The conditions on $P\/$ imply that $P\/$ admits a $w\/$-labeling where the points of $P-\{w \}\/$ are $v_0, v_1, \ldots, v_{2n-2}\/$ (in anticlockwise direction).  Let $T_0\/$ be a $w\/$-caterpillar  on $P\/$ constructed in the following way.

\vspace{3mm} {\em Case (1):} $T_0\/$ is Type-$1\/$.

\vspace{1mm}
(i)  Let  $ w v_{n-1} v_0 v_{n}\/$ be a $4\/$-path of   $T_0\/$.

\vspace{1mm}
(ii) Let $v_s, v_t\/$ be two points of $P - \{w \}\/$ such that the line passing through $v_s,v_t\/$ is the $r\/$-halving line. Here  $n+1 \leq s < t \leq 2n-1\/$ (where $v_{2n-1} = v_0\/$). Then either choose $v_s x\/$ to be an edge of $T_0\/$ for every $x \in \{v_{s+1}, v_{s+2}, \ldots, v_{t}\} \/$,  or else choose $v_t y\/$ to be an edge of $T_0\/$ for every  $ y \in \{v_{s}, v_{s+1}, \ldots, v_{t-1} \}\/$.

\vspace{1mm}
(iii) With respect to the edge $v_sv_t\/$ of $T_0\/$, we either choose $v_tv_{s-1}\/$ or else $v_sv_{t+1}\/$ where $v_{t+1} \not \in \{v_1,v_2, \ldots, v_{n-1}\} \cup \{w\}\/$ to be an edge of $T_0\/$. Iteratively, if $v_pv_q\/$ has been chosen to be an edge of $T_0\/$, where $0 \leq p \leq s\/$, $t \leq q \leq n-1\/$, then we either choose $v_qv_{p-1}\/$ or else $v_pv_{q+1}\/$ (where $v_{q+1} \not \in \{v_1,v_2, \ldots, v_{n-1}\} \cup \{w\}\/$) to be an edge of $T_0\/$.

\vspace{1mm}
(iv) Now to  every edge $v_iv_j \/$ that has been chosen to be in  $T_0\/$ where $i, j \in \{n+1, n+2, \ldots, 2n-1\}\/$, choose $v_{i-n}v_{j-n}\/$ to be an edge of $T_0\/$.

\vspace{1mm}
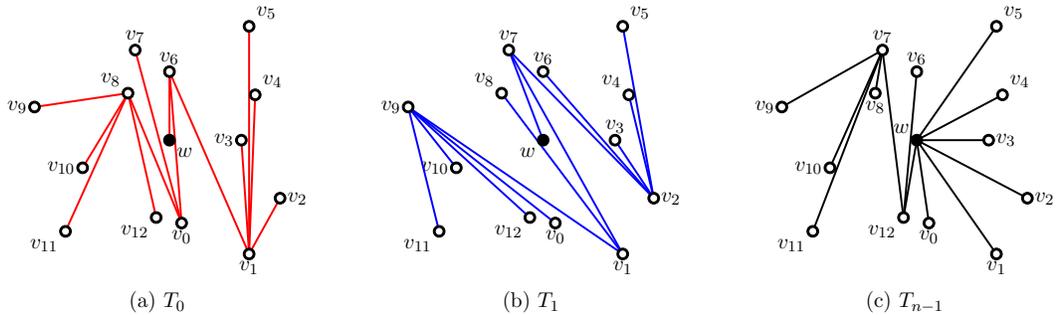
\begin{figure}[htb]
\resizebox{15cm}{!}{
\begin{minipage}{.45\textwidth}
\begin{tikzpicture}

\coordinate (center) at (0,0);
   \foreach \x in {0,27.692307,...,360} {

               }
       \coordinate (v13) at (0,0);\node [below right] at (v13) {${w}$};
       \coordinate (v0) at (0.22,-1.5);\node[below] at (v0) {${v_0}$};
       \coordinate (v1) at (1.44,-2.06);\node[below] at (v1) {${v_1}$};
       \coordinate (v2) at (2,-1.05);\node[ right] at (v2) {${v_2}$};
       \coordinate (v3) at (1.3,0);\node[left] at (v3) {${v_3}$};
       \coordinate (v4) at (1.55,0.82);\node[above right] at (v4) {${v_4}$};
       \coordinate (v5) at (1.44,2.06);\node[above right] at (v5) {${v_5}$};
       \coordinate (v6) at (0,1.24);\node[above] at (v6) {${v_6}$};
       \coordinate (v7) at (-0.62,1.63);\node[above] at (v7) {${v_7}$};
       \coordinate (v8) at (-0.75,0.85);\node[above left] at (v8) {${v_{8}}$};
       \coordinate (v9) at (-2.44,0.6);\node[left] at (v9) {${v_{9}}$};
       \coordinate (v10) at (-1.57,-0.5);\node[left] at (v10) {${v_{10}}$};
       \coordinate (v11) at (-1.88,-1.65);\node[below left] at (v11) {${v_{11}}$};
       \coordinate (v12) at (-0.24,-1.4);\node[below left] at (v12) {${v_{12}}$};

       \draw [line width=1,red](v6) -- (v13);
       \draw [line width=1,red](v6) -- (v0);
       \draw [line width=1,red](v0) -- (v7);
        \draw [line width=1,red](v0) -- (v8);
         \draw [line width=1,red](v8) -- (v9); \draw [line width=1,red](v8) -- (v10);\draw [line width=1,red](v8) -- (v11);\draw [line width=1,red](v8) -- (v12);
        \draw [line width=1,red](v6) -- (v1);
        \draw [line width=1,red](v1) -- (v2);\draw [line width=1,red](v1) -- (v3);\draw [line width=1,red](v1) -- (v4);\draw [line width=1,red](v1) -- (v5);


       \filldraw[black] (v13) circle(3pt);\filldraw[black] (v0) circle(3pt);\filldraw[black] (v1) circle(3pt);\filldraw[black] (v2) circle(3pt);
       \filldraw[black] (v3) circle(3pt);\filldraw[black] (v4) circle(3pt);\filldraw[black] (v5) circle(3pt);\filldraw[black] (v6) circle(3pt);
       \filldraw[black] (v7) circle(3pt);\filldraw[black] (v8) circle(3pt);\filldraw[black] (v9) circle(3pt);
       \filldraw[black] (v10) circle(3pt);\filldraw[black] (v11) circle(3pt);\filldraw[black] (v12) circle(3pt);

       \filldraw[white] (v0) circle(1.5pt);\filldraw[white] (v1) circle(1.5pt);\filldraw[white] (v2) circle(1.5pt);
       \filldraw[white] (v3) circle(1.5pt);\filldraw[white] (v4) circle(1.5pt);\filldraw[white] (v5) circle(1.5pt);\filldraw[white] (v6) circle(1.5pt);
       \filldraw[white] (v7) circle(1.5pt);\filldraw[white] (v8) circle(1.5pt);\filldraw[white] (v9) circle(1.5pt);
       \filldraw[white] (v10) circle(1.5pt);\filldraw[white] (v11) circle(1.5pt);\filldraw[white] (v12) circle(1.5pt);

\end{tikzpicture}
\centering

(a) $T_0$
\end{minipage}
\begin{minipage}{.45\textwidth}
\begin{tikzpicture}
\coordinate (center) at (0,0);
   \foreach \x in {0,27.692307,...,360} {

               }
       \coordinate (v13) at (0,0);\node [below left] at (v13) {${w}$};
       \coordinate (v0) at (0.22,-1.5);\node[below] at (v0) {${v_0}$};
       \coordinate (v1) at (1.44,-2.06);\node[below] at (v1) {${v_1}$};
       \coordinate (v2) at (2,-1.05);\node[ right] at (v2) {${v_2}$};
       \coordinate (v3) at (1.3,0);\node[above] at (v3) {${v_3}$};
       \coordinate (v4) at (1.55,0.82);\node[above left] at (v4) {${v_4}$};
       \coordinate (v5) at (1.44,2.06);\node[above right] at (v5) {${v_5}$};
       \coordinate (v6) at (0,1.24);\node[above] at (v6) {${v_6}$};
       \coordinate (v7) at (-0.62,1.63);\node[above] at (v7) {${v_7}$};
       \coordinate (v8) at (-0.75,0.85);\node[above left] at (v8) {${v_{8}}$};
       \coordinate (v9) at (-2.44,0.6);\node[left] at (v9) {${v_{9}}$};
       \coordinate (v10) at (-1.57,-0.5);\node[left] at (v10) {${v_{10}}$};
       \coordinate (v11) at (-1.88,-1.65);\node[below left] at (v11) {${v_{11}}$};
       \coordinate (v12) at (-0.24,-1.4);\node[below left] at (v12) {${v_{12}}$};


       \draw [line width=1,blue](v1) -- (v7);
       \draw [line width=1,blue](v1) -- (v8);
       \draw [line width=1,blue](v7) -- (v13);
       \draw [line width=1,blue](v1) -- (v9);
       \draw [line width=1,blue](v9) -- (v10);\draw [line width=1,blue](v9) -- (v11);\draw [line width=1,blue](v9) -- (v12);\draw [line width=1,blue](v9) -- (v0);
       \draw [line width=1,blue](v7) -- (v2);
       \draw [line width=1,blue](v2) -- (v3);\draw [line width=1,blue](v2) -- (v4);\draw [line width=1,blue](v2) -- (v5);\draw [line width=1,blue](v2) -- (v6);


       \filldraw[black] (v13) circle(3pt);\filldraw[black] (v0) circle(3pt);\filldraw[black] (v1) circle(3pt);\filldraw[black] (v2) circle(3pt);
       \filldraw[black] (v3) circle(3pt);\filldraw[black] (v4) circle(3pt);\filldraw[black] (v5) circle(3pt);\filldraw[black] (v6) circle(3pt);
       \filldraw[black] (v7) circle(3pt);\filldraw[black] (v8) circle(3pt);\filldraw[black] (v9) circle(3pt);
       \filldraw[black] (v10) circle(3pt);\filldraw[black] (v11) circle(3pt);\filldraw[black] (v12) circle(3pt);

       \filldraw[white] (v0) circle(1.5pt);\filldraw[white] (v1) circle(1.5pt);\filldraw[white] (v2) circle(1.5pt);
       \filldraw[white] (v3) circle(1.5pt);\filldraw[white] (v4) circle(1.5pt);\filldraw[white] (v5) circle(1.5pt);\filldraw[white] (v6) circle(1.5pt);
       \filldraw[white] (v7) circle(1.5pt);\filldraw[white] (v8) circle(1.5pt);\filldraw[white] (v9) circle(1.5pt);
       \filldraw[white] (v10) circle(1.5pt);\filldraw[white] (v11) circle(1.5pt);\filldraw[white] (v12) circle(1.5pt);

\end{tikzpicture}
\centering

(b) $T_1$
\end{minipage}
\begin{minipage}{.45\textwidth}
\begin{tikzpicture}
\coordinate (center) at (0,0);
   \foreach \x in {0,27.692307,...,360} {

               }
       \coordinate (v13) at (0,0);\node [above left] at (v13) {${w}$};
       \coordinate (v0) at (0.22,-1.5);\node[below] at (v0) {${v_0}$};
       \coordinate (v1) at (1.44,-2.06);\node[below] at (v1) {${v_1}$};
       \coordinate (v2) at (2,-1.05);\node[ right] at (v2) {${v_2}$};
       \coordinate (v3) at (1.3,0);\node[right] at (v3) {${v_3}$};
       \coordinate (v4) at (1.55,0.82);\node[above right] at (v4) {${v_4}$};
       \coordinate (v5) at (1.44,2.06);\node[above right] at (v5) {${v_5}$};
       \coordinate (v6) at (0,1.24);\node[above] at (v6) {${v_6}$};
       \coordinate (v7) at (-0.62,1.63);\node[above] at (v7) {${v_7}$};
       \coordinate (v8) at (-0.75,0.85);\node[below] at (v8) {${v_{8}}$};
       \coordinate (v9) at (-2.44,0.6);\node[left] at (v9) {${v_{9}}$};
       \coordinate (v10) at (-1.57,-0.5);\node[left] at (v10) {${v_{10}}$};
       \coordinate (v11) at (-1.88,-1.65);\node[below left] at (v11) {${v_{11}}$};
       \coordinate (v12) at (-0.24,-1.4);\node[below left] at (v12) {${v_{12}}$};

      \draw [line width=1,black](v12) -- (v13);\draw [line width=1,black](v0) -- (v13);
      \draw [line width=1,black](v1) -- (v13);\draw [line width=1,black](v2) -- (v13);\draw [line width=1,black](v3) -- (v13);
      \draw [line width=1,black](v4) -- (v13);\draw [line width=1,black](v5) -- (v13);
      \draw [line width=1,black](v6) -- (v12);

      \draw [line width=1,black](v12) -- (v7);
      \draw [line width=1,black](v7) -- (v8);\draw [line width=1,black](v7) -- (v9);\draw [line width=1,black](v7) -- (v10);\draw [line width=1,black](v7) -- (v11);

       \filldraw[black] (v13) circle(3pt);\filldraw[black] (v0) circle(3pt);\filldraw[black] (v1) circle(3pt);\filldraw[black] (v2) circle(3pt);
       \filldraw[black] (v3) circle(3pt);\filldraw[black] (v4) circle(3pt);\filldraw[black] (v5) circle(3pt);\filldraw[black] (v6) circle(3pt);
       \filldraw[black] (v7) circle(3pt);\filldraw[black] (v8) circle(3pt);\filldraw[black] (v9) circle(3pt);
       \filldraw[black] (v10) circle(3pt);\filldraw[black] (v11) circle(3pt);\filldraw[black] (v12) circle(3pt);

       \filldraw[white] (v0) circle(1.5pt);\filldraw[white] (v1) circle(1.5pt);\filldraw[white] (v2) circle(1.5pt);
       \filldraw[white] (v3) circle(1.5pt);\filldraw[white] (v4) circle(1.5pt);\filldraw[white] (v5) circle(1.5pt);\filldraw[white] (v6) circle(1.5pt);
       \filldraw[white] (v7) circle(1.5pt);\filldraw[white] (v8) circle(1.5pt);\filldraw[white] (v9) circle(1.5pt);
       \filldraw[white] (v10) circle(1.5pt);\filldraw[white] (v11) circle(1.5pt);\filldraw[white] (v12) circle(1.5pt);

\end{tikzpicture}
\centering

(c) $T_{n-1}\/$
\end{minipage}
}
\caption{ An example of Type-$1$ $w\/$-caterpillar $T_0\/$  and its derived trees in Theorem \ref{the4}.   \label{wheel-1}  }
\end{figure}

\vspace{3mm} {\em Case (2):} $T_0\/$ is Type-2.

\vspace{1mm}
(i) Let  $ w v_{n-1} v_0 v_{n}\/$ be a $4\/$-path of   $T_0\/$.

\vspace{1mm}
(ii) Let $v_s, v_t\/$ be two points of $P - \{w \}\/$ such that the line passing through $v_s,v_t\/$ is the $r\/$-halving line. Here  $n \leq s < t \leq 2n-2\/$. Then either choose $v_s x\/$ to be an edge of $T_0\/$ for every $x \in \{v_{s+1}, v_{s+2}, \ldots, v_{t}\} \/$,  or else choose $v_t y\/$ to be an edge of $T_0\/$ for every  $ y \in \{v_{s}, v_{s+1}, \ldots, v_{t-1} \}\/$.

\vspace{1mm}
(iii) With respect to the edge $v_sv_t\/$ of $T_0\/$, we either choose $v_tv_{s-1}\/$ (where $v_{s-1} \not \in \{v_0,v_1, \ldots, v_{n-1}\} \cup \{w\}\/$) or else $v_sv_{t+1}\/$ to be an edge of $T_0\/$. Iteratively, if $v_pv_q\/$ has been chosen to be an edge of $T_0\/$, where $n \leq p \leq s\/$, $t \leq q \leq 2n-2\/$, then we either choose $v_qv_{p-1}\/$ (where $v_{p-1} \not \in \{v_0,v_1, \ldots, v_{n-1}\} \cup \{w\}\/$) or else $v_pv_{q+1}\/$ to be an edge of $T_0\/$.

\vspace{1mm}
(iv) Now to  every edge $v_iv_j \/$ that has been chosen to be in  $T_0\/$ where $i, j \in \{n, n+1, \ldots, 2n-2\}\/$, choose $v_{i-n}v_{j-n}\/$ to be an edge of $T_0\/$.

\vspace{1mm} It can be shown in the way similar to the case of Theorem \ref{the3} each of the above  $T_0\/$ is a plane tree.

\vspace{1mm}
For each $q \in \{1, 2, \ldots, n-2\}\/$, let $T_q\/$ denote the tree on $P\/$ defined by $v_iv_j \in E(T_q)\/$ if and only if $v_{i-q}v_{j-q} \in E(T_0)\/$. Clearly $T_q\/$ is isomorphic to $T_0\/$. The proof that  $E(T_0) \cap E(T_q) = \emptyset\/$ is similar to that of Theorem \ref{the3}.


\vspace{1mm} Finally, let $T_{n-1}\/$ be the caterpillar on $P\/$ constructed in the following way.

\begin{figure}[htb]
\resizebox{15cm}{!}{
\begin{minipage}{.45\textwidth}
\begin{tikzpicture}
\coordinate (center) at (0,0);
   \foreach \x in {0,27.692307,...,360} {

               }
       \coordinate (v13) at (0,0);\node [below right] at (v13) {${w}$};
       \coordinate (v0) at (0.22,-1.5);\node[below] at (v0) {${v_0}$};
       \coordinate (v1) at (1.44,-2.06);\node[below] at (v1) {${v_1}$};
       \coordinate (v2) at (2,-1.05);\node[ right] at (v2) {${v_2}$};
       \coordinate (v3) at (1.3,0);\node[below left] at (v3) {${v_3}$};
       \coordinate (v4) at (1.55,0.82);\node[above right] at (v4) {${v_4}$};
       \coordinate (v5) at (1.44,2.06);\node[above right] at (v5) {${v_5}$};
       \coordinate (v6) at (0,1.24);\node[above] at (v6) {${v_6}$};
       \coordinate (v7) at (-0.62,1.63);\node[above] at (v7) {${v_7}$};
       \coordinate (v8) at (-0.75,0.85);\node[above left] at (v8) {${v_{8}}$};
       \coordinate (v9) at (-2.44,0.6);\node[left] at (v9) {${v_{9}}$};
       \coordinate (v10) at (-1.57,-0.5);\node[left] at (v10) {${v_{10}}$};
       \coordinate (v11) at (-1.88,-1.65);\node[below left] at (v11) {${v_{11}}$};
       \coordinate (v12) at (-0.24,-1.4);\node[below left] at (v12) {${v_{12}}$};

       \draw [line width=1,red](v6) -- (v13);
       \draw [line width=1,red](v6) -- (v0);
       \draw [line width=1,red](v0) -- (v7);
       \draw [line width=1,red](v7) -- (v12);\draw [line width=1,red](v0) -- (v5);
       \draw [line width=1,red](v12) -- (v8);\draw [line width=1,red](v12) -- (v9);\draw [line width=1,red](v12) -- (v10);   \draw [line width=1,red](v12) -- (v11);
       \draw [line width=1,red](v5) -- (v1);\draw [line width=1,red](v5) -- (v2);\draw [line width=1,red](v5) -- (v3);\draw [line width=1,red](v5) -- (v4);


       \filldraw[black] (v13) circle(3pt);\filldraw[black] (v0) circle(3pt);\filldraw[black] (v1) circle(3pt);\filldraw[black] (v2) circle(3pt);
       \filldraw[black] (v3) circle(3pt);\filldraw[black] (v4) circle(3pt);\filldraw[black] (v5) circle(3pt);\filldraw[black] (v6) circle(3pt);
       \filldraw[black] (v7) circle(3pt);\filldraw[black] (v8) circle(3pt);\filldraw[black] (v9) circle(3pt);
       \filldraw[black] (v10) circle(3pt);\filldraw[black] (v11) circle(3pt);\filldraw[black] (v12) circle(3pt);

       \filldraw[white] (v0) circle(1.5pt);\filldraw[white] (v1) circle(1.5pt);\filldraw[white] (v2) circle(1.5pt);
       \filldraw[white] (v3) circle(1.5pt);\filldraw[white] (v4) circle(1.5pt);\filldraw[white] (v5) circle(1.5pt);\filldraw[white] (v6) circle(1.5pt);
       \filldraw[white] (v7) circle(1.5pt);\filldraw[white] (v8) circle(1.5pt);\filldraw[white] (v9) circle(1.5pt);
       \filldraw[white] (v10) circle(1.5pt);\filldraw[white] (v11) circle(1.5pt);\filldraw[white] (v12) circle(1.5pt);

\end{tikzpicture}
\centering

(a) $T_0$
\end{minipage}
\begin{minipage}{.45\textwidth}
\begin{tikzpicture}
\coordinate (center) at (0,0);
   \foreach \x in {0,27.692307,...,360} {

               }
       \coordinate (v13) at (0,0);\node [below left] at (0.1,0) {${w}$};
       \coordinate (v0) at (0.22,-1.5);\node[below right] at (v0) {${v_0}$};
       \coordinate (v1) at (1.44,-2.06);\node[below] at (v1) {${v_1}$};
       \coordinate (v2) at (2,-1.05);\node[ right] at (v2) {${v_2}$};
       \coordinate (v3) at (1.3,0);\node[right] at (v3) {${v_3}$};
       \coordinate (v4) at (1.55,0.82);\node[above right] at (v4) {${v_4}$};
       \coordinate (v5) at (1.44,2.06);\node[above right] at (v5) {${v_5}$};
       \coordinate (v6) at (0,1.24);\node[above] at (v6) {${v_6}$};
       \coordinate (v7) at (-0.62,1.63);\node[above] at (v7) {${v_7}$};
       \coordinate (v8) at (-0.75,0.85);\node[above left] at (v8) {${v_{8}}$};
       \coordinate (v9) at (-2.44,0.6);\node[left] at (v9) {${v_{9}}$};
       \coordinate (v10) at (-1.57,-0.5);\node[left] at (v10) {${v_{10}}$};
       \coordinate (v11) at (-1.88,-1.65);\node[below left] at (v11) {${v_{11}}$};
       \coordinate (v12) at (-0.24,-1.4);\node[left] at (v12) {${v_{12}}$};


       \draw [line width=1,blue](v1) -- (v7);
       \draw [line width=1,blue](v1) -- (v8);
       \draw [line width=1,blue](v7) -- (v13);
      \draw [line width=1,blue](v8) -- (v0);
      \draw [line width=1,blue](v0) -- (v9); \draw [line width=1,blue](v0) -- (v10); \draw [line width=1,blue](v0) -- (v11); \draw [line width=1,blue](v0) -- (v12);
       \draw [line width=1,blue](v1) -- (v6);
        \draw [line width=1,blue](v6) -- (v2);\draw [line width=1,blue](v6) -- (v3);\draw [line width=1,blue](v6) -- (v4);\draw [line width=1,blue](v6) -- (v5);

       \filldraw[black] (v13) circle(3pt);\filldraw[black] (v0) circle(3pt);\filldraw[black] (v1) circle(3pt);\filldraw[black] (v2) circle(3pt);
       \filldraw[black] (v3) circle(3pt);\filldraw[black] (v4) circle(3pt);\filldraw[black] (v5) circle(3pt);\filldraw[black] (v6) circle(3pt);
       \filldraw[black] (v7) circle(3pt);\filldraw[black] (v8) circle(3pt);\filldraw[black] (v9) circle(3pt);
       \filldraw[black] (v10) circle(3pt);\filldraw[black] (v11) circle(3pt);\filldraw[black] (v12) circle(3pt);

       \filldraw[white] (v0) circle(1.5pt);\filldraw[white] (v1) circle(1.5pt);\filldraw[white] (v2) circle(1.5pt);
       \filldraw[white] (v3) circle(1.5pt);\filldraw[white] (v4) circle(1.5pt);\filldraw[white] (v5) circle(1.5pt);\filldraw[white] (v6) circle(1.5pt);
       \filldraw[white] (v7) circle(1.5pt);\filldraw[white] (v8) circle(1.5pt);\filldraw[white] (v9) circle(1.5pt);
       \filldraw[white] (v10) circle(1.5pt);\filldraw[white] (v11) circle(1.5pt);\filldraw[white] (v12) circle(1.5pt);

\end{tikzpicture}
\centering

(b) $T_1$
\end{minipage}
\begin{minipage}{.45\textwidth}
\begin{tikzpicture}
\coordinate (center) at (0,0);
   \foreach \x in {0,27.692307,...,360} {

               }
      \coordinate (v13) at (0,0);\node [above left] at (v13) {${w}$};
       \coordinate (v0) at (0.22,-1.5);\node[below] at (v0) {${v_0}$};
       \coordinate (v1) at (1.44,-2.06);\node[below] at (v1) {${v_1}$};
       \coordinate (v2) at (2,-1.05);\node[ right] at (v2) {${v_2}$};
       \coordinate (v3) at (1.3,0);\node[right] at (v3) {${v_3}$};
       \coordinate (v4) at (1.55,0.82);\node[above right] at (v4) {${v_4}$};
       \coordinate (v5) at (1.44,2.06);\node[above right] at (v5) {${v_5}$};
       \coordinate (v6) at (0,1.24);\node[above] at (v6) {${v_6}$};
       \coordinate (v7) at (-0.62,1.63);\node[above] at (v7) {${v_7}$};
       \coordinate (v8) at (-0.75,0.85);\node[above left] at (v8) {${v_{8}}$};
       \coordinate (v9) at (-2.44,0.6);\node[left] at (v9) {${v_{9}}$};
       \coordinate (v10) at (-1.57,-0.5);\node[above left] at (v10) {${v_{10}}$};
       \coordinate (v11) at (-1.88,-1.65);\node[below left] at (v11) {${v_{11}}$};
       \coordinate (v12) at (-0.24,-1.4);\node[below left] at (v12) {${v_{12}}$};

      \draw [line width=1,black](v12) -- (v13);\draw [line width=1,black](v0) -- (v13);
      \draw [line width=1,black](v1) -- (v13);\draw [line width=1,black](v2) -- (v13);\draw [line width=1,black](v3) -- (v13);
      \draw [line width=1,black](v4) -- (v13);\draw [line width=1,black](v5) -- (v13);
      \draw [line width=1,black](v6) -- (v12);

      \draw [line width=1,black](v6) -- (v11);
     \draw [line width=1,black](v11) -- (v7);\draw [line width=1,black](v11) -- (v8);\draw [line width=1,black](v11) -- (v9);\draw [line width=1,black](v11) -- (v10);

       \filldraw[black] (v13) circle(3pt);\filldraw[black] (v0) circle(3pt);\filldraw[black] (v1) circle(3pt);\filldraw[black] (v2) circle(3pt);
       \filldraw[black] (v3) circle(3pt);\filldraw[black] (v4) circle(3pt);\filldraw[black] (v5) circle(3pt);\filldraw[black] (v6) circle(3pt);
       \filldraw[black] (v7) circle(3pt);\filldraw[black] (v8) circle(3pt);\filldraw[black] (v9) circle(3pt);
       \filldraw[black] (v10) circle(3pt);\filldraw[black] (v11) circle(3pt);\filldraw[black] (v12) circle(3pt);

       \filldraw[white] (v0) circle(1.5pt);\filldraw[white] (v1) circle(1.5pt);\filldraw[white] (v2) circle(1.5pt);
       \filldraw[white] (v3) circle(1.5pt);\filldraw[white] (v4) circle(1.5pt);\filldraw[white] (v5) circle(1.5pt);\filldraw[white] (v6) circle(1.5pt);
       \filldraw[white] (v7) circle(1.5pt);\filldraw[white] (v8) circle(1.5pt);\filldraw[white] (v9) circle(1.5pt);
       \filldraw[white] (v10) circle(1.5pt);\filldraw[white] (v11) circle(1.5pt);\filldraw[white] (v12) circle(1.5pt);

\end{tikzpicture}
\centering

(c) $T_{n-1}\/$
\end{minipage}
}
\caption{ An example of Type-$2$ $w\/$-caterpillar $T_0\/$  and its derived trees in Theorem \ref{the4}.   \label{wheel-2}  }
\end{figure}
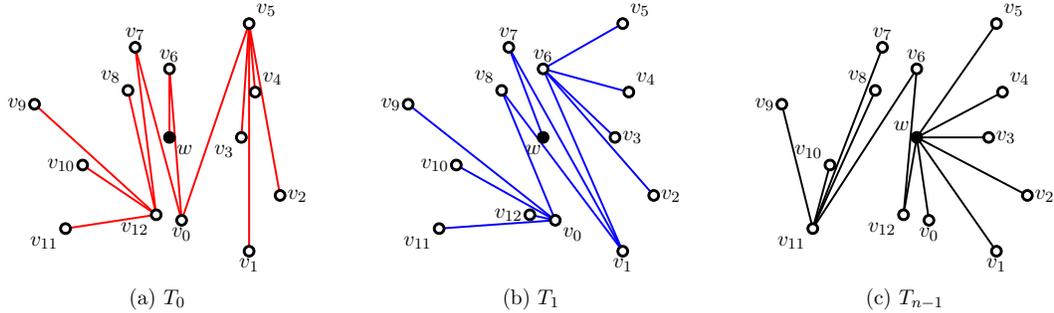

\vspace{1mm} Let $v_{i+n-1}v_{j+n-1} \in E(T_{n-1})\/$ if and only if $v_iv_j \in E(T_0)\/$ where $i,j \in \{0,1, \ldots, n-1\}\/$ and $wz \in E(T_{n-1})\/$ for every $z \in \{v_{2n-2}, v_0, \ldots, v_{n-2}\} \/$ (see for example the tree $T_{n-1}\/$ in Figures \ref{wheel-1} and  \ref{wheel-2}). Clearly, $T_{n-1}\/$ is a plane spanning tree on $P\/$. \qed

\vspace{1mm}

\end{document}